\documentclass[11pt,a4paper]{article}
\usepackage[OT1]{fontenc}
\usepackage{amsmath,amssymb,amsthm}
\usepackage{graphicx}
\usepackage{booktabs,array}
\usepackage[margin=25mm]{geometry}
\usepackage{microtype}
\usepackage{enumitem}
\usepackage{xcolor}
\usepackage[colorlinks=true,linkcolor=blue!50!black,citecolor=blue!50!black,urlcolor=blue!50!black]{hyperref}
\hypersetup{pdftitle={Unilluminable rooms: minimal examples, classifications and exceptions},
  pdfauthor={J. C. F. de Winter and George W. Tokarsky}}

\theoremstyle{plain}
\newtheorem{theorem}{Theorem}[section]
\newtheorem{lemma}[theorem]{Lemma}
\newtheorem{proposition}[theorem]{Proposition}
\newtheorem{corollary}[theorem]{Corollary}
\theoremstyle{definition}
\newtheorem{definition}[theorem]{Definition}
\newtheorem{example}[theorem]{Example}
\theoremstyle{remark}
\newtheorem{remark}[theorem]{Remark}

\newcommand{\Z}{\mathbb{Z}}
\newcommand{\T}{\mathcal{T}}
\newcommand{\Sq}{\mathcal{S}}
\newcommand{\K}{\mathcal{K}}
\newcommand{\E}{\mathcal{E}}
\newcommand{\dg}{^\circ}
\setlist{itemsep=2pt,topsep=4pt}
\numberwithin{equation}{section}
\title{Unilluminable rooms:\\[2pt] minimal examples, classifications and exceptions}
\author{J. C. F. de Winter\\[2pt]
  {\small Department of Cognitive Robotics}\\
  {\small Delft University of Technology}\\
  {\small Delft, The Netherlands}\\
  {\small\href{mailto:j.c.f.dewinter@tudelft.nl}{j.c.f.dewinter@tudelft.nl}}
  \and George W. Tokarsky\\[2pt]
  {\small Department of Mathematical and Statistical Sciences}\\
  {\small University of Alberta}\\
  {\small Edmonton, Alberta, Canada}\\
  {\small\href{mailto:tokarsky@ualberta.ca}{tokarsky@ualberta.ca}}}
\date{31 October 2026}

\begin{document}
\maketitle

\begin{abstract}
A room with mirrored walls can contain two points that are mutually invisible: no light ray from one reaches the other. Tokarsky built the first polygonal examples by tiling a room with reflected copies of a triangle and making every tile vertex of certain kinds a corner of the room. We study the smallest rooms of this kind, which we call lattice rooms. Their interior tile vertices of the remaining kind, the sources, are invisible and pairwise mutually invisible. On the tiling by isosceles right triangles, a hole-free lattice room with a source has at least 12 sides, and the pinwheel is the only one with 12 sides; on the square tiling the least number is 16, again attained by a single room. A hole-free lattice room with two sources has at least 24 sides on the first tiling and at least 34 on the square tiling; up to congruence there are 12 such rooms with 24 sides, 89 with 25 sides and 14 with 34 sides. A single absorbing post, a single partition, or a single hole that touches the outer boundary at one point lowers the least number of sides for two sources to 18, 20 and 22 (21 if sides are counted as maximal straight segments); two and three posts give 17 and 16 sides, and a hole that does not touch the outer boundary gives 26. We list all smallest rooms in each case. Rooms assembled from $n$ modules have exactly $12n$ sides and $n$ sources, and a hole-free lattice room with three or more sources has at least 36 sides. We also give rooms in which some tile vertices need not be corners: 27-sided rooms on the tiling by triangles with angles $30\dg$, $60\dg$ and $90\dg$, which are smallest in their class and form 15 types, a nonagon with three mutually invisible corners, and a 42-sided room in which a rotation through $180\dg$ takes the place of a corner. Finally, in a ring of three trapezoids a light ray leaves a corner and returns to it, although that corner is invisible in a single trapezoid. Several classifications rely on exhaustive computer searches made finite by the bounds proved here; programs and exact coordinates accompany the paper.
\end{abstract}

\section{Introduction}\label{sec:intro}

Imagine a room whose walls are mirrors, lit by a single candle. Does the light reach every point of the room? The question is attributed to Ernst Straus in the early 1950s~\cite{Tok95}, and Klee~\cite{Klee} first published it, for polygonal rooms. L. S. Penrose and R. Penrose~\cite{Penrose} described a room with curved walls that cannot be lit completely from any point. For polygonal rooms the first examples are due to Tokarsky~\cite{Tok95}: a room with 26 sides that contains two points $P$ and $Q$ such that no light ray from $P$ reaches $Q$. We call such points \emph{mutually invisible}, and a room that contains them \emph{unilluminable}. Stewart~\cite{Stewart} described this room in \emph{Scientific American}. Castro~\cite{Castro}, Tokarsky~\cite[Figure~2]{Tok96} and Thomson~\cite{Thomson} then found such rooms with 24 sides; Tokarsky's room appeared in a later note by Stewart~\cite{Stewart97}. Rooms with 22 sides exist if a hole may touch the outer boundary at one point: Tokarsky~\cite[Figure~3(b)]{Tok96} drew one with a triangular hole, Thomson~\cite{Thomson} found the same room and one with a square hole, which has 21 sides when sides are counted as maximal straight segments, and Wolecki~\cite{Wolecki} published the room with the triangular hole. Wolecki also reports an unpublished 22-sided example of Monteil.

Such pairs are uncommon. Leli\`evre, Monteil and Weiss~\cite{LMW} proved that in a polygonal room in which the angle between any two sides is a rational multiple of $\pi$, each point is mutually invisible with only finitely many points, and Wolecki~\cite{Wolecki} proved that such a room contains only finitely many mutually invisible pairs. All rooms in this paper are of this kind. In Tokarsky's construction~\cite{Tok95} such pairs come from the arithmetic of a tiling. The room is a union of isosceles right triangles taken from the tiling of the plane generated by reflections in their sides; $P$ and $Q$ are lattice points of one kind; and every lattice point of the other kinds in the room is a corner of the room. A light ray from $P$ to $Q$ unfolds into a straight segment between two lattice points, this segment passes through a lattice point of another kind, and that point folds back onto a corner of the room, where the ray is absorbed.

This paper asks how few sides such a room can have, and how many rooms attain the least number. In a letter to Stewart~\cite{Tok96}, the second author drew several small examples: a pinwheel with 12 sides and a cross with 16 sides, each with an invisible point at its centre, and rooms with two mutually invisible points that use a post, a partition or a hole. We show that these examples are smallest in natural classes, and we list every room of the smallest size.

We work with \emph{lattice rooms}: unions of tiles of a reflection tiling in which every tile vertex outside a distinguished lattice is a corner (Definition~\ref{def:latticeroom}). The points of that lattice in the interior of the room are its \emph{sources}, and they are invisible and pairwise mutually invisible. Two counting lemmas (Lemmas~\ref{lem:euler} and~\ref{lem:wings}) relate the numbers of tiles, sides and sources of a lattice room, so that each classification becomes a finite problem. Table~\ref{tab:summary} summarises the results.

\begin{table}[htbp]
\centering\small
\begin{tabular}{@{}>{\raggedright\arraybackslash}p{0.36\textwidth}>{\raggedright\arraybackslash}p{0.12\textwidth}c>{\centering\arraybackslash}p{0.12\textwidth}l@{}}
\toprule
Rooms & Sources & Least $N$ & Rooms with least $N$ & Result\\
\midrule
Hole-free lattice rooms on $\T$ & one & 12 & 1 & Theorem~\ref{thm:pinwheel}\\
Hole-free lattice rooms on $\Sq$ & one & 16 & 1 & Theorem~\ref{thm:square}\\
Rooms on $\T$ with one post & two & 18 & 5 & Theorem~\ref{thm:post}\\
Rooms on $\T$ with two posts & two & 17 & 6 & Theorem~\ref{thm:posts}\\
Rooms on $\T$ with three posts & two & 16 & 4 & Theorem~\ref{thm:posts}\\
Rooms on $\T$ with one partition & two & 20 & 8 & Corollary~\ref{cor:partition}\\
Rooms on $\T$ with one pinched hole & two & 22 & 4 & Theorem~\ref{thm:pinched}\\
Lattice rooms on $\T$ with one hole & two & 26 & 12 & Theorem~\ref{thm:hole}\\
Hole-free lattice rooms on $\T$ & two & 24 & 12 & Theorem~\ref{thm:tri}\\
Hole-free lattice rooms on $\Sq$ & two & 34 & 14 & Theorem~\ref{thm:sq}\\
Hole-free lattice rooms on $\T$ & three & 36 & & Theorem~\ref{thm:three}\\
Hole-free rooms on $\K$ whose $C$-points are corners & two & 27 & 15 & Theorem~\ref{thm:min27}\\
\bottomrule
\end{tabular}
\caption{Least number $N$ of sides in each class, and the number of rooms with $N$ sides, counted up to congruence. Here $\T$, $\Sq$ and $\K$ are the tilings by isosceles right triangles, by squares, and by triangles with angles $30\dg$, $60\dg$, $90\dg$. In addition, there are 21 post rooms with 19 sides and 89 hole-free lattice rooms on $\T$ with 25 sides. The rooms with three sources and 36 sides were not counted. In a pinched room the point where the hole meets the outer boundary is an endpoint of sides of both polygons (Section~\ref{sec:devices}); counted as maximal straight segments, the least number is 21 (Remark~\ref{rem:pinchcount}).}
\label{tab:summary}
\end{table}

Section~\ref{sec:many} assembles larger rooms from modules: an assembly of $n$ modules has exactly $12n$ sides and $n$ sources, which are invisible and pairwise mutually invisible, and a hole-free lattice room with three sources has at least $36$ sides. Section~\ref{sec:relax} relaxes the corner condition. The blocking argument uses only part of it, and we give the $15$ rooms with $27$ sides on the tiling by triangles with angles $30\dg$, $60\dg$ and $90\dg$, which are smallest in their class, a nonagon with three mutually invisible corners of angle $90\dg$, and a 42-sided room in which a rotation through $180\dg$ takes the place of a corner. Section~\ref{sec:trapezoid} shows how a hole can exchange the labels of a tile: in a ring of three trapezoids, a light ray leaves a corner and returns to it, although this corner is invisible in a single trapezoid.

The classifications in Sections~\ref{sec:pinwheel}, \ref{sec:devices}, \ref{sec:holefree} and~\ref{sec:27} and the bound of Theorem~\ref{thm:three} use exhaustive computer searches. In each case the bounds proved in the text make the search finite, and Section~\ref{sec:computations} describes the computations. The programs and the exact coordinates of all rooms accompany the paper.

\section{Light rays, tilings and lattice rooms}\label{sec:prelim}

\subsection{Rooms and light rays}
A \emph{room} is a compact set $R$ in the plane with connected interior whose boundary is a finite union of segments. Except in Section~\ref{sec:devices}, the boundary of a room consists of finitely many pairwise disjoint simple closed polygons: the outer boundary, and the boundaries of the \emph{holes}. A \emph{corner} of $R$ is a boundary point near which the boundary is not a straight segment, and a \emph{side} is a maximal straight segment of the boundary. When the boundary consists of disjoint simple polygons, the number $N$ of sides equals the number of corners. Section~\ref{sec:devices} states how sides are counted in rooms with posts, a partition or a pinched hole.

The boundary of a room is a mirror. A \emph{light ray} is a polygonal path of positive length in $R$ that is straight between reflections and reflects at the boundary with equal angles of incidence and reflection. It reflects only at points in the relative interiors of sides; a ray that reaches a corner is absorbed there and stops, as in~\cite{Klee,Tok95}. A ray may start at any point of $R$, including a boundary point or a corner. Two points $P$ and $Q$ of $R$ are \emph{mutually invisible} if no light ray starts at $P$ and ends at $Q$; by reversibility this relation is symmetric. Following Tokarsky and Boese~\cite{TB}, we call a point $P$ \emph{invisible} if no light ray starts and ends at $P$.

\subsection{Tilings}
Let $\T$ be the tiling of the plane in which each unit square $[i,i+1]\times[j,j+1]$, with $i,j\in\Z$, is cut into two isosceles right triangles by its diagonal through $(i,j)$ if $i+j$ is even, and by its other diagonal if $i+j$ is odd. The edges of $\T$ lie on the \emph{mirror lines}
\begin{equation}\label{eq:mirrors}
x=k,\qquad y=k,\qquad x-y=2k,\qquad x+y=2k\qquad(k\in\Z),
\end{equation}
and every tile is the mirror image of each neighbour in their common edge. We call a lattice point $(x,y)\in\Z^2$ an \emph{$A$-point} if $x$ and $y$ are both even, a \emph{$B$-point} if both are odd, and a \emph{$C$-point} otherwise. Every tile of $\T$ has one vertex of each kind, with its right angle at the $C$-vertex. Eight tiles meet at each $A$- and $B$-point, and four at each $C$-point. Two lattice points are \emph{adjacent} if their distance is $1$; the four $C$-points adjacent to an $A$-point $(x,y)$ are $(x\pm1,y)$ and $(x,y\pm1)$.

Let $\Sq$ be the tiling of the plane by the unit squares $[i,i+1]\times[j,j+1]$; its mirror lines are the lines $x=k$ and $y=k$. Tokarsky~\cite[Example~4]{Tok95} used this tiling for rooms whose sides are all horizontal or vertical. We call points of $\Z^2$ $A$-, $B$- and $C$-points as before. Every square has one $A$-vertex, one $B$-vertex and two $C$-vertices.

\begin{definition}\label{def:latticeroom}
A \emph{lattice room} on $\T$ (or on $\Sq$) is a room $R$ that is a union of tiles of $\T$ (or of $\Sq$), whose boundary consists of pairwise disjoint simple closed polygons, and which satisfies the \emph{corner condition}: every $B$-point and every $C$-point of $R$ is a corner of $R$. The \emph{sources} of $R$ are the $A$-points in its interior. A lattice point on the boundary that is not a corner is \emph{flat}.
\end{definition}

By the corner condition, every interior lattice point of a lattice room is a source and every flat point is an $A$-point. We write $T$ for the number of tiles, $M$ for the number of sources, $f$ for the number of flat points and $h$ for the number of holes. We count rooms up to congruence: two rooms have the same \emph{type} if an isometry of the plane maps one onto the other (for the rooms of Section~\ref{sec:devices}, together with the posts, the partition or the point where the hole meets the outer boundary).

\subsection{The blocking lemma}
All our examples rest on an argument of Tokarsky~\cite[Lemma~3.1 and Theorem~3.2]{Tok95}; see also~\cite[Corollary]{TB}. On $\T$ it runs as follows. Let $P$ and $Q$ be $A$-points of a lattice room. Unfolding a light ray from $P$ to $Q$ turns it into a straight segment from $P$ to another $A$-point. Such a segment passes through a $B$- or $C$-point, so the light ray passes through a $B$- or $C$-point of the room. That point is a corner, where the ray stops. We state the argument for a general lattice, since Section~\ref{sec:relax} applies it to other tilings.

A \emph{lattice} $\Lambda$ is the set of points $me_1+ne_2$ with $m,n\in\Z$, for two linearly independent vectors $e_1$ and $e_2$. We write $2\Lambda$ for the set of points $me_1+ne_2$ with $m$ and $n$ both even, and $\Lambda\setminus2\Lambda$ for the other points of $\Lambda$, those with $m$ or $n$ odd. For $\Lambda=\Z^2$, with $e_1=(1,0)$ and $e_2=(0,1)$, the points of $2\Lambda$ are the $A$-points, and the points of $\Lambda\setminus2\Lambda$ are the $B$- and $C$-points. Only some reflections map $\Lambda$ onto $\Lambda$ and $2\Lambda$ onto $2\Lambda$. For example, the reflection $(x,y)\mapsto(2k-x,y)$ in the line $x=k$ does this for $\Lambda=\Z^2$ when $k$ is an integer; for $k=\frac12$ it maps $\Z^2$ onto itself and the $A$-point $(0,0)$ to the $C$-point $(1,0)$. Condition~(i) of the next lemma asks for this property of the line of every side of the room, and condition~(ii) asks for corners; both are hypotheses, which Corollary~\ref{cor:blocking} verifies for lattice rooms.

\begin{lemma}[Blocking lemma]\label{lem:blocking}
Let $\Lambda$ be a lattice, and let $R$ be a room with the following two properties.
\begin{enumerate}[label=\textup{(\roman*)}]
\item For every side of $R$, the reflection in the line of that side maps $\Lambda$ onto $\Lambda$ and $2\Lambda$ onto $2\Lambda$.
\item Every point of $\Lambda\setminus2\Lambda$ that lies in $R$ is a corner of $R$ or a point at which light is absorbed.
\end{enumerate}
Then every point of $2\Lambda$ that lies in $R$ is invisible, and any two such points are mutually invisible.
\end{lemma}

\begin{proof}
Suppose that a light ray $\gamma$ starts at a point $P$ of $2\Lambda$ and ends at a point $Q$ of $2\Lambda$, where $P=Q$ is allowed. We unfold $\gamma$. At the first point where $\gamma$ reflects, reflect the remaining part of $\gamma$ in the line of that side, so that the first two pieces of the path lie on one straight line. Repeat this at each later reflection point of the new path. The result is a straight segment $[P,Q']$ of the same length as $\gamma$. A point $X$ of $[P,Q']$ corresponds to a point $\phi(X)$ of $\gamma$, where $\phi$ undoes the reflections made before the piece that contains $X$; in particular $\phi(Q')=Q$. Each such $\phi$ is a product of reflections in lines of sides of $R$, so by~(i) it maps $\Lambda$ onto $\Lambda$ and $2\Lambda$ onto $2\Lambda$, and hence also $\Lambda\setminus2\Lambda$ onto $\Lambda\setminus2\Lambda$. In particular $Q'$ lies in $2\Lambda$.

Since $P$ and $Q'$ lie in $2\Lambda$, we can write $Q'-P=2ae_1+2be_2$ with integers $a$ and $b$, which are not both $0$ because $\gamma$ has positive length. Let $g$ be the greatest common divisor of $a$ and $b$, and let $X=P+(a/g)e_1+(b/g)e_2$. Then $X$ lies on $[P,Q']$ at the fraction $1/(2g)\le\frac12$ of the way from $P$ to $Q'$, so strictly between them. The point $X$ lies in $\Lambda$, and it does not lie in $2\Lambda$, because $P$ does and the integers $a/g$ and $b/g$ are not both even. So $\phi(X)$ is a point of $\Lambda\setminus2\Lambda$ through which $\gamma$ passes before it ends. By~(ii) it is a corner of $R$ or a point that absorbs light, and a light ray stops at such a point. This contradiction proves the lemma.
\end{proof}

\begin{corollary}\label{cor:blocking}
In a lattice room on $\T$ or on $\Sq$, every $A$-point is invisible, and any two $A$-points are mutually invisible.
\end{corollary}

\begin{proof}
Take $\Lambda=\Z^2$, so that $2\Lambda$ is the set of $A$-points and $\Lambda\setminus2\Lambda$ is the set of $B$- and $C$-points. Condition~(ii) of Lemma~\ref{lem:blocking} is then the corner condition. For condition~(i), every side of a lattice room lies on a mirror line~\eqref{eq:mirrors}, and the reflections in these lines are
\[
(x,y)\mapsto(2k-x,\,y),\quad(x,\,2k-y),\quad(y+2k,\,x-2k),\quad(2k-y,\,2k-x).
\]
Each of them keeps the parities of the two coordinates or interchanges them, so it maps $\Z^2$ onto $\Z^2$ and the $A$-points onto the $A$-points. On $\Sq$ only the first two kinds of lines occur.
\end{proof}

\subsection{Counting tiles}

\begin{lemma}\label{lem:euler}
Let $R$ be a union of $T$ tiles of $\T$ with connected interior whose boundary consists of $h+1$ pairwise disjoint simple closed polygons, and let $I$ and $b$ be the numbers of lattice points in the interior and on the boundary of $R$. Then $b=T-2I+2-2h$. For a union of $T$ squares of $\Sq$ the same is true with $T$ replaced by $2T$. In particular, a lattice room satisfies
\begin{equation}\label{eq:euler}
N=T-2M+2-2h-f\quad\text{on }\T,\qquad N=2T-2M+2-2h-f\quad\text{on }\Sq.
\end{equation}
\end{lemma}

\begin{proof}
The vertices and edges of the tiles of $R$ form a connected plane graph with $I+b$ vertices, $E$ edges and $T+h+1$ faces: the $T$ tiles, the $h$ holes and the outer face. The boundary polygons have $b$ edges in all, so counting pairs of a tile and one of its edges gives $3T=2E-b$ on $\T$ and $4T=2E-b$ on $\Sq$. Euler's formula $(I+b)-E+(T+h+1)=2$, that is, $(I+b)-E+T=1-h$, and substituting $E$ prove the first statement. In a lattice room $I=M$ and $b=N+f$.
\end{proof}

\begin{lemma}\label{lem:wings}
Let $R$ be a lattice room on $\T$ or on $\Sq$.
\begin{enumerate}[label=\textup{(\alph*)}]
\item Each $C$-point of $R$ lies in exactly one or exactly three tiles of $R$.
\item Let $s$ be a source. Each of the four $C$-points adjacent to $s$ lies in exactly three tiles of $R$, and exactly one of these three tiles does not contain $s$. We call these four tiles the \emph{wings} of $s$; together with the tiles at $s$, eight on $\T$ and four on $\Sq$, they form the \emph{module} of $s$.
\item $R$ satisfies $T\ge12M+3f$ on $\T$, and $2T\ge12M+3f$ on $\Sq$.
\end{enumerate}
\end{lemma}

\begin{proof}
(a) Four tiles meet at a $C$-point $c$. If $R$ contains exactly two of them, these are adjacent and $c$ is flat, or opposite and the boundary of $R$ is not simple at $c$; if $R$ contains all four, $c$ is interior. Each case contradicts Definition~\ref{def:latticeroom}.

(b) The tiles at $s$ supply two tiles at each adjacent $C$-point, so by (a) there is exactly one more.

(c) We argue on $\T$; on $\Sq$ replace $T$ by $2T$, since each square has two $C$-vertices. The four tiles at a $C$-point $c$ form two pairs, the tiles of a pair sharing the edge from $c$ to one of the two $A$-points adjacent to $c$. Call such an edge \emph{full} if both tiles containing it lie in $R$, and let $U$ be the number of full edges. Every tile has exactly one $C$-vertex, so $T$ is the sum, over all $C$-points, of the number of tiles of $R$ there. By (a) this number is $0$, $1$ or $3$. At a $C$-point with three tiles exactly one of its two edges is full, and at a $C$-point with one tile none is. Every full edge has exactly one endpoint that is a $C$-point, so $U$ is the number of $C$-points with three tiles, and $T=3U+Q$, where $Q$ is the number of $C$-points in exactly one tile of $R$. In particular $T\ge3U$. At a source all four edges to adjacent $C$-points are full, and at a flat point at least one is, since the tiles of $R$ at a flat point fill a half-plane. Every full edge has exactly one endpoint that is an $A$-point, so these edges are distinct, and $U\ge4M+f$. We write $U=4M+f+\delta$ with $\delta\ge0$. Hence $T\ge3U=12M+3f+3\delta\ge12M+3f$.
\end{proof}

\begin{corollary}\label{cor:bound}
A hole-free lattice room on $\T$ or on $\Sq$ satisfies $N\ge10M+2+2f$. In particular it has at least $12$ sides if it has a source, and at least $22$ sides if it has two.
\end{corollary}

\begin{proof}
Combine~\eqref{eq:euler}, with $h=0$, and Lemma~\ref{lem:wings}(c).
\end{proof}

The next lemma shows that congruent rooms on our tilings are related by symmetries of the lattice.

\begin{lemma}\label{lem:congruence}
Let $R$ and $R'$ be nonempty finite unions of tiles of $\T$, or of $\Sq$. Then every isometry $g$ with $g(R)=R'$ maps $\Z^2$ onto itself, and so maps the interior lattice points of $R$ onto those of $R'$.
\end{lemma}

\begin{proof}
The boundary of $R$ is a union of tile edges. A maximal straight segment in it is therefore horizontal or vertical with integer length, or diagonal with length an integer multiple of $\sqrt2$, and its endpoints are lattice points. The same is true for $R'$, and $g$ maps these segments of $R$ onto those of $R'$ and preserves their lengths. Since $\sqrt2$ is irrational, the linear part of $g$ maps the direction of an axis to that of an axis, or the direction of a diagonal to that of a diagonal. Being orthogonal, it is then one of the eight symmetries of a square centred at the origin, and it maps $\Z^2$ onto itself. As $g$ also maps an endpoint of a segment to a lattice point, $g$ maps $\Z^2$ onto itself.
\end{proof}

We use Lemma~\ref{lem:congruence} to count types. In each classification below, the listed rooms have their sources at the same points. By Lemma~\ref{lem:congruence}, a congruence $g$ between two listed rooms maps interior lattice points onto interior lattice points. In a lattice room these are the sources. In the rooms of Section~\ref{sec:devices} they can also include the posts and the free end of a partition, and $g$ maps these onto the posts and the free end. So $g$ maps the set of sources onto itself. As $g$ also maps $\Z^2$ onto itself and maps a source to a source, it maps $A$-points onto $A$-points, and so it is a symmetry of the tiling. Two listed rooms therefore have the same type exactly when a symmetry of the tiling that maps the set of sources onto itself maps one room onto the other. These symmetries form a finite group, and the types are its orbits: the classes of listed rooms that the group maps onto one another.

The searches of Sections~\ref{sec:devices}--\ref{sec:many} are made finite by a bound on the distance between sources. The \emph{sup-distance} of two points $(x,y)$ and $(x',y')$ is $\max(|x-x'|,|y-y'|)$. The tiles at a lattice point $s$ lie in the square of side $2$ centred at $s$, and the module of a source $s$ lies in the square of side $4$ centred at $s$.

\begin{lemma}\label{lem:chain}
Let $r\ge1$ be an integer, and let $t_0,t_1,\dots,t_{j+1}$ be tiles of $\T$, or of $\Sq$, such that $t_i$ and $t_{i+1}$ have a common edge for $0\le i\le j$. If $t_0$ lies in the square of side $2r$ centred at a lattice point $s$, and $t_{j+1}$ lies in the square of side $2r$ centred at a lattice point $s'$, then $s$ and $s'$ have sup-distance at most $j+2r$.
\end{lemma}

\begin{proof}
Each tile lies in a unit square $[a,a+1]\times[b,b+1]$ with $a,b\in\Z$. Two tiles with a common edge lie in the same unit square, if the edge is a diagonal, or in two unit squares with a common side, so the values of $a$ for consecutive tiles $t_i$ and $t_{i+1}$ differ by at most $1$. Let $a_0$ and $a_{j+1}$ be the values for $t_0$ and $t_{j+1}$. Then $|a_{j+1}-a_0|\le j+1$, and $s_x-r\le a_0\le s_x+r-1$ and $s'_x-r\le a_{j+1}\le s'_x+r-1$, so $|s'_x-s_x|\le|a_{j+1}-a_0|+2r-1\le j+2r$. The same argument applies to the $y$-coordinates.
\end{proof}

If a room has connected interior, any two disjoint sets of its tiles are joined by a chain of its tiles as in Lemma~\ref{lem:chain}, and in a shortest such chain the inner tiles $t_1,\dots,t_j$ lie outside both sets.

\section{One source}\label{sec:one}

\subsection{The pinwheel}\label{sec:pinwheel}
\begin{theorem}\label{thm:pinwheel}
Every hole-free lattice room on $\T$ with a source has at least $12$ sides. Up to congruence, the only such room with $12$ sides is the \emph{pinwheel} of Figure~\textup{\ref{fig:one}(a)}, with vertices
\[
(-2,0),\ (-1,0),\ (-1,1),\ (0,2),\ (0,1),\ (1,1),\ (2,0),\ (1,0),\ (1,-1),\ (0,-2),\ (0,-1),\ (-1,-1).
\]
Its source $(0,0)$ is invisible.
\end{theorem}

\begin{proof}
By Corollary~\ref{cor:bound}, $N\ge10M+2+2f\ge12$, and $N=12$ forces $M=1$ and $f=0$, hence $T=12$ by~\eqref{eq:euler}. Place the source at the origin. The room contains the eight tiles at the origin and, by Lemma~\ref{lem:wings}(b), the four wings of the source; these twelve tiles are all of $R$. The wing at $(1,0)$ is one of the triangles with vertices $(1,0)$, $(2,0)$, $(1,\pm1)$, so it contains exactly one of the $B$-points $(1,1)$ and $(1,-1)$; we say that the wing points to this $B$-point, and similarly for the other wings. Each of the four $B$-points $(\pm1,\pm1)$ lies in two tiles at the origin and is adjacent to two of the four $C$-points $(\pm1,0)$, $(0,\pm1)$. If both wings at these $C$-points point to the same $B$-point, then four consecutive tiles meet there and it is flat, contrary to the corner condition. So each $B$-point receives at most one wing, hence exactly one. Going around the origin, the wings therefore all point clockwise or all point counterclockwise. The two rooms obtained are mirror images of each other, and each is the pinwheel. Conversely, the pinwheel satisfies the corner condition and has 12 sides, and Corollary~\ref{cor:blocking} shows that its source is invisible.
\end{proof}

The corners of the pinwheel are its four $A$-points $(\pm2,0)$, $(0,\pm2)$, its four $B$-points and its four $C$-points; its area is $6$, and it has rotational symmetry of order four about its source. The pinwheel is Figure~1(b) of the second author's letter~\cite{Tok96} and room~a of~\cite{Stewart97}, where it is conjectured that no room with fewer sides contains an invisible point; Theorem~\ref{thm:pinwheel} confirms this for hole-free lattice rooms on $\T$. Tokarsky and Boese~\cite[Figure~6]{TB} note that the source and the four $A$-corners of the pinwheel are pairwise mutually invisible.

Larger lattice rooms on $\T$ with one source need not contain the pinwheel. In the 13-gon with vertices
\[
\begin{gathered}(-2,-1),\ (-1,-1),\ (0,-2),\ (0,-1),\ (1,-1),\ (2,0),\ (1,0),\\ (1,1),\ (0,1),\ (0,2),\ (-1,1),\ (-1,0),\ (-2,0),\end{gathered}
\]
the source is $(0,0)$, the wings at $(-1,0)$ and $(0,-1)$ both point to the $B$-point $(-1,-1)$, and a thirteenth tile, with vertices $(-2,-1)$, $(-2,0)$ and $(-1,-1)$, makes this point a corner. By Corollary~\ref{cor:bound} and~\eqref{eq:euler}, a hole-free lattice room on $\T$ with one source and $N$ sides has $f\le(N-12)/2$ flat points and $T=N+f$ tiles. An exhaustive search over such tile sets shows that up to congruence there are exactly three of these rooms with $13$ sides (Figure~\ref{fig:onesource13}), and none of them contains the pinwheel. For $N=12,13,\dots,20$ the numbers of types are $1$, $3$, $4$, $19$, $67$, $84$, $253$, $783$ and $1346$.

\subsection{The cross}
\begin{theorem}\label{thm:square}
Every hole-free lattice room on $\Sq$ with a source has at least $16$ sides. Up to congruence, the only such room with $16$ sides is the \emph{cross}
\[
\Sigma=[-1,1]^2\cup[-2,-1]\times[0,1]\cup[-1,0]\times[1,2]\cup[1,2]\times[-1,0]\cup[0,1]\times[-2,-1]
\]
of Figure~\textup{\ref{fig:one}(b)}.
\end{theorem}

\begin{proof}
Let $R$ be such a room with $N\le16$. By~\eqref{eq:euler}, $f=2T-2M+2-N$, and Corollary~\ref{cor:bound} gives $2f\le N-10M-2$; eliminating $f$ yields $4T\le3N-6M-6$. As $N\ge10M+2$, we have $M=1$, and then $T\le9$.

Place the source at the origin. Its four squares form $[-1,1]^2$, and $R$ contains the four wings, one beyond each $C$-point $(\pm1,0)$, $(0,\pm1)$. Each wing contains exactly one of the $B$-points $(\pm1,\pm1)$, and each $B$-point lies in one square of $[-1,1]^2$ and is adjacent to two of the four $C$-points. If a $B$-point lies in exactly one wing, then two adjacent squares of $R$ meet there, and $R$ needs a further square at that point to make it a corner. Each wing contains one $B$-point and each $B$-point lies in at most two wings, so the numbers of wings at the four $B$-points are $0$, $1$ or $2$ and add up to $4$; hence the number of $B$-points that lie in exactly one wing is even. If it is nonzero, $R$ needs at least two further squares, since a unit square cannot contain two of the points $(\pm1,\pm1)$; then $T\ge10$, which is impossible. Hence each $B$-point lies in zero or two wings. So two $B$-points lie in two wings each and the other two lie in none. Adjacent $B$-points, such as $(1,1)$ and $(1,-1)$, cannot both lie in two wings, since the single wing at the $C$-point between them would have to contain both. So the two $B$-points are opposite, which up to a quarter turn leaves one possibility. The eight squares obtained form $\Sigma$.

Suppose that $R$ has a ninth square. It shares an edge with $\Sigma$, since the interior of $R$ is connected, and one endpoint of this edge is a $C$-point. This $C$-point lies in one or three squares of $\Sigma$ and would lie in two or four squares of $R$, contrary to Lemma~\ref{lem:wings}(a). Thus $R=\Sigma$, which has $16$ sides.
\end{proof}

The corners of $\Sigma$ are four $A$-points, four $B$-points and eight $C$-points, and its area is $8$. Its symmetry group has order four and is generated by the reflections in the two diagonals. The cross is Figure~1(a) of the second author's letter~\cite{Tok96}, and Tokarsky and Boese~\cite[Figure~6]{TB} note that its five $A$-points are pairwise mutually invisible.

\begin{proposition}\label{prop:sigma}
In $\Sigma$, the five $A$-points $(0,0)$, $(\pm2,0)$, $(0,\pm2)$ are invisible and pairwise mutually invisible. The source $(0,0)$ is the only invisible interior point of $\Sigma$.
\end{proposition}

\begin{proof}
The first statement is Corollary~\ref{cor:blocking}. Let $X=(x,y)$ be an interior point of $\Sigma$. If $x\notin\Z$, the vertical line through $X$ leaves $\Sigma$ through horizontal sides, at points that are not corners because every corner has integer coordinates. A ray from $X$ along this line meets such a side at a right angle and returns to $X$. If $x\in\Z$ and $y\notin\Z$, a horizontal ray does the same. The only interior lattice point of $\Sigma$ is the origin.
\end{proof}

Every other interior point of $\Sigma$ is joined to the origin by a straight light ray, since $\Sigma$ is star-shaped about the origin and the open segment from the origin to an interior point contains no corner.

\subsection{A star on another tiling}
\begin{figure}[htbp]
\centering
\includegraphics[width=\textwidth]{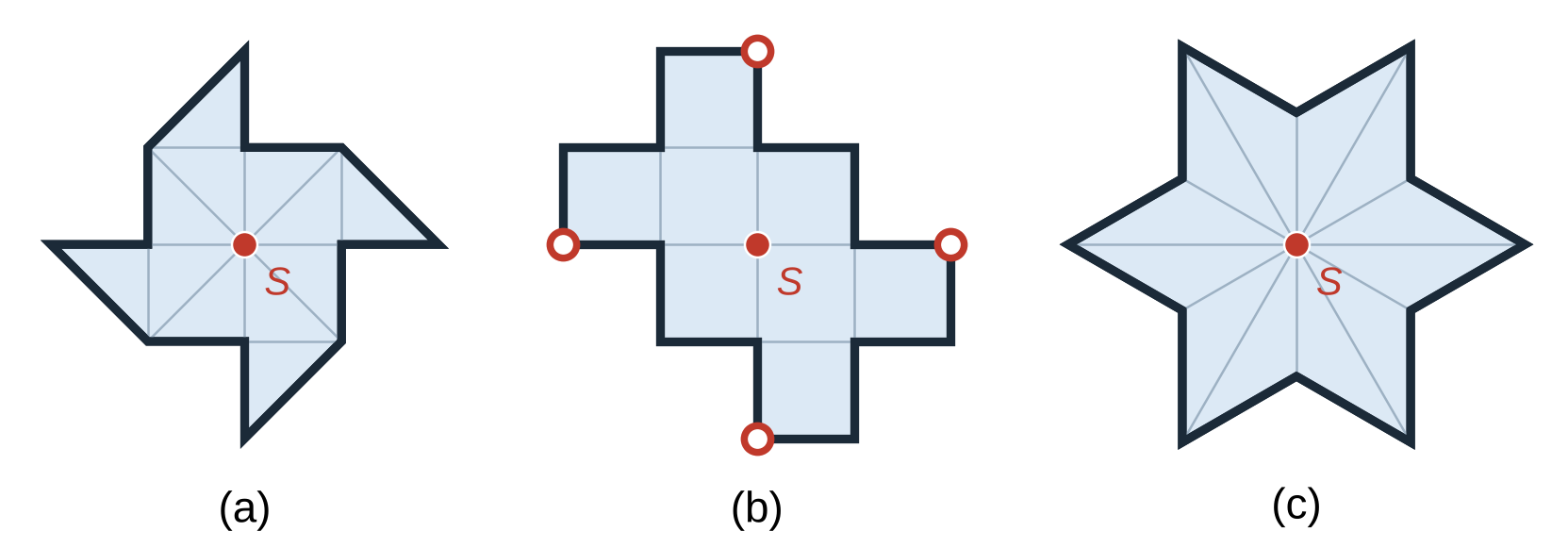}
\caption{Rooms with an invisible interior point $S$. (a)~The pinwheel, on $\T$, with 12 sides. (b)~The cross $\Sigma$, on $\Sq$, with 16 sides; the four open circles mark the $A$-points on its boundary, which together with $S$ are pairwise mutually invisible. (c)~The star on the tiling $\E$, with 12 sides. Thin lines are tile edges.}
\label{fig:one}
\end{figure}

The same idea works on other tilings, where the corner condition can leave only one room. Let $\E$ be the tiling obtained from the regular tiling by equilateral triangles by dividing each triangle into three triangles with angles $30\dg$, $30\dg$, $120\dg$, using the segments from its centre to its vertices. Fix a vertex $S$ of the equilateral tiling. The \emph{star} of $S$ is the union of the twelve tiles of $\E$ that have $S$ as a vertex (Figure~\ref{fig:one}(c)). It has twelve sides, and its corners alternate between six vertices of the equilateral tiling, with angle $60\dg$, and six centres of equilateral triangles, with angle $240\dg$.

\begin{proposition}\label{prop:star}
\begin{enumerate}[label=\textup{(\alph*)}]
\item The point $S$ is invisible in the star.
\item Let $R$ be a union of tiles of $\E$ with connected interior that contains the star of $S$, and suppose that every centre of an equilateral triangle that lies in $R$ is a corner of $R$. Then $R$ is the star of $S$.
\end{enumerate}
\end{proposition}

\begin{proof}
(a) Any two tiles of the star with a common edge are mirror images of each other in that edge, and $S$ is a vertex with angle $30\dg$ of every tile. Fix one tile $\Delta$. \emph{Folding} maps each tile onto $\Delta$ by the reflections in the edges crossed along a chain of tiles from that tile to $\Delta$; going once around $S$ crosses twelve edges at angles of $30\dg$, whose reflections compose to the identity, so the result does not depend on the chain. Folding maps a light ray in the star that meets $S$ at most at its ends to a light ray in $\Delta$, and it maps the star onto the triangle $\Delta$ with angles $\pi/6$ at a vertex $A$, $\pi/6$ at a vertex $B$ and $2\pi/3$ at a vertex $C$, so that $S$ goes to $A$, the other vertices of the equilateral tiling go to $B$, and the centres go to $C$. A light ray in the star from $S$ back to $S$, cut at its first return to $S$, folds into a light ray in $\Delta$ from $A$ back to $A$; it avoids $B$ and $C$, because their preimages in the star are corners. By a lemma of Tokarsky~\cite{Tok95a},~\cite[Lemma~4.1]{Tok95}, which applies because the angle $\pi/6$ at $A$ divides $\pi/2$ and the angle at $B$ is a multiple of it, no light ray in $\Delta$ starts and ends at $A$; see also~\cite[Theorem~3]{TB} and~\cite[Lemma~3.1]{Wolecki}.

(b) Each side of the star joins a vertex of the equilateral tiling to the centre $z$ of an equilateral triangle that has $S$ as a vertex, and the tile of $\E$ beyond this side is the third tile of that triangle at $z$. If $R$ contained it, $z$ would be an interior point of $R$. Since the interior of $R$ is connected, any tile of $R$ outside the star would lead, through a chain of tiles sharing edges, to such a tile. Hence $R$ is the star.
\end{proof}

Thus the star admits no wings: on $\E$, the corner condition at the centres leaves the bare star as the only room of this kind containing it.

\section{Two sources with posts, a partition or a hole}\label{sec:devices}

The second author's letter~\cite[Figure~3]{Tok96} describes three ways to separate two sources with fewer sides: an absorbing post at a single point, a thin partition attached to the outer boundary, and a hole. In this section we find the smallest rooms of each kind on $\T$ and list them all. We also treat rooms with several posts and rooms with a hole that does not touch the outer boundary.

\begin{definition}\label{def:devices}
Let $R$ be a union of tiles of $\T$ with connected interior.
\begin{enumerate}[label=\textup{(\alph*)}]
\item $R$ is a \emph{post room} if its boundary is a simple closed polygon and there is a point $p$ in its interior, the \emph{post}, such that every $B$- and $C$-point of $R$ other than $p$ is a corner. Light that reaches $p$ is absorbed.
\item $R$ is a \emph{partition room} if its boundary is a simple closed polygon and $R$ contains a segment $J$, the \emph{partition}, with the following properties: $J$ lies on a mirror line; one endpoint of $J$ is a corner of $R$; the other endpoint $p$, the \emph{free end}, and the relative interior of $J$ lie in the interior of $R$; and the relative interior of $J$ contains no $A$-point. Every $B$- and $C$-point of $R$ other than $p$ must be a corner of $R$. Both faces of $J$ are mirrors, and light that reaches $p$ is absorbed.
\item $R$ is a \emph{pinched room} if its boundary consists of two simple closed polygons, the outer boundary and the boundary of the hole, that meet in exactly one point $q$, and every $B$- and $C$-point of $R$ other than $q$ is a corner. Light that reaches $q$ is absorbed.
\end{enumerate}
The \emph{sources} of $R$ are the $A$-points in its interior other than the post or the free end.
\end{definition}

We count sides as follows. A post adds no sides, and a partition adds two, one for each face. In a pinched room we count the sides of the outer polygon and of the hole polygon separately, and in each polygon we treat $q$ as an endpoint of sides, even when that polygon is straight at $q$. Remark~\ref{rem:pinchcount} explains this convention.

Lemma~\ref{lem:blocking} applies with $\Lambda=\Z^2$, as in the proof of Corollary~\ref{cor:blocking}. Every side lies on a mirror line, including each face of a partition, which gives condition~(i); and every $B$- and $C$-point of $R$ is a corner or absorbs light, which gives condition~(ii). So the sources of a post room, a partition room or a pinched room are invisible and pairwise mutually invisible.

\subsection{One post}
\begin{figure}[htbp]
\centering
\includegraphics[width=0.78\textwidth]{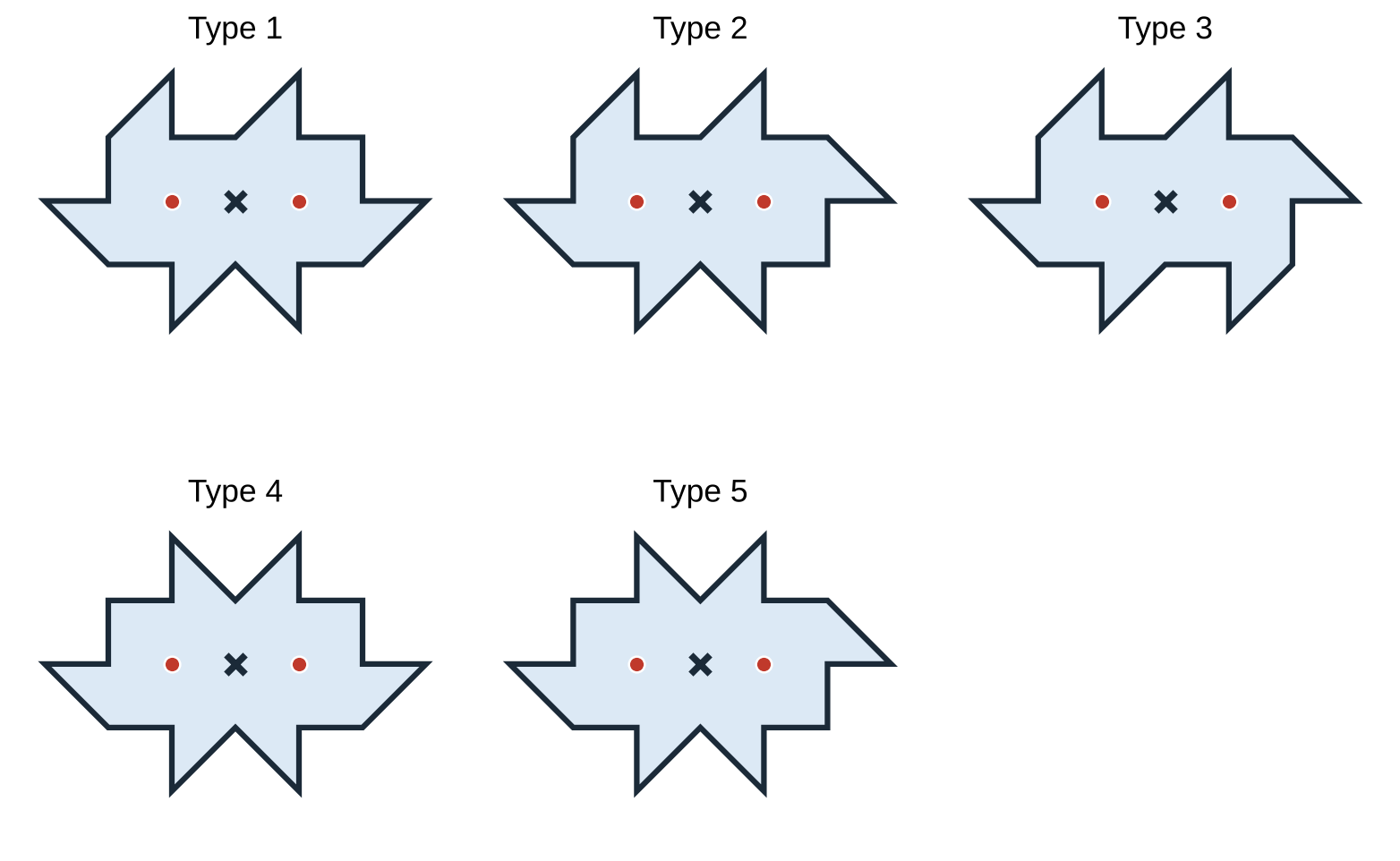}
\caption{The five post rooms with 18 sides (Theorem~\ref{thm:post}). The cross marks the post $(1,0)$; the sources are $(0,0)$ and $(2,0)$. Type~3 was found by the second author.}
\label{fig:post18}
\end{figure}

\begin{theorem}\label{thm:post}
A post room with at least two sources has at least $18$ sides. Up to congruence there are exactly five post rooms with two sources and $18$ sides, and exactly twenty-one with $19$ sides \textup{(}Figures~\textup{\ref{fig:post18}} and~\textup{\ref{fig:post19}}\textup{)}. In each of them, after a symmetry of $\T$, the sources are $(0,0)$ and $(2,0)$ and the post is $(1,0)$.
\end{theorem}

\begin{proof}
Let $R$ be a post room with $M\ge2$ sources and $N\le19$ sides. If the post $p$ is not a $B$- or $C$-point, then $R$ is a lattice room and $N\ge22$ by Corollary~\ref{cor:bound}. So $p$ is a $B$- or $C$-point, the interior lattice points of $R$ are the $M$ sources and $p$, and every flat point is an $A$-point. Lemma~\ref{lem:euler} gives
\begin{equation}\label{eq:posteuler}
N+f=T-2(M+1)+2.
\end{equation}
As in the proof of Lemma~\ref{lem:wings}(c), $T$ is the sum over $C$-points of the number of tiles of $R$ there. This number is $0$, $1$ or $3$, except that it is $4$ at $p$ if $p$ is a $C$-point; then both edges at $p$ are full. Let $Q$ be the number of $C$-points in exactly one tile, let $\varepsilon=1$ if $p$ is a $C$-point and $\varepsilon=0$ otherwise, and let $U$ be the number of full edges; as in Lemma~\ref{lem:wings}(c), $U=4M+f+\delta$ with $\delta\ge0$. Then $T=3U+Q-2\varepsilon$, and~\eqref{eq:posteuler} becomes
\[
N=10M+2f+Q+3\delta-2\varepsilon.
\]
If $p$ is a $B$-point, $N\ge20$. Hence $p$ is a $C$-point, and $N=10M+2f+Q+3\delta-2\le19$ forces $M=2$ and $f=0$; moreover $N\ge18$, and $T=N+4\le23$.

The sixteen tiles at the two sources are distinct. Suppose that the sources have no adjacent $C$-point in common. Then each of their eight adjacent $C$-points lies in two of these tiles and needs at least one more, and these further tiles are distinct because their $C$-vertices differ. This gives $T\ge24$, a contradiction. Hence the sources share an adjacent $C$-point. After a symmetry of $\T$ they are $(0,0)$ and $(2,0)$, and $(1,0)$ lies in four tiles of $R$. It is therefore an interior point, and it is the post.

The other six $C$-points adjacent to the sources, $(-1,0)$, $(0,\pm1)$, $(2,\pm1)$ and $(3,0)$, each lie in two tiles at the sources and need exactly one wing. For $N=18$ we have $T=22$, so $R$ consists of the sixteen tiles at the sources and six wings, one chosen out of two at each of these points: $2^6=64$ candidates. For $N=19$ we have $T=23$, and $R$ contains one further tile, which shares an edge with the other 22 tiles because the interior of $R$ is connected; this gives $576$ distinct candidates. We tested each candidate, with exact integer arithmetic, for a simple boundary, for the corner condition away from $(1,0)$ and for the number of sides. This leaves $14$ rooms with $18$ sides and $84$ with $19$ sides. The isometries that map $\{(0,0),(2,0)\}$ onto itself are the identity, the reflections in the lines $y=0$ and $x=1$, and the rotation through $180\dg$ about $(1,0)$; they are symmetries of $\T$ and form a group of order four. By Lemma~\ref{lem:congruence}, a congruence between two of our rooms maps interior lattice points to interior lattice points; as it also maps the post to the post, it maps sources to sources and so lies in this group. The types are the orbits of this group, and Table~\ref{tab:post} gives their number and size.
\end{proof}

\begin{table}[htbp]
\centering\small
\begin{tabular}{@{}ccccl@{}}
\toprule
Sides & Candidates & Rooms & Types & Rooms per type\\
\midrule
18 & 64 & 14 & 5 & 2, 2, 2, 4, 4\\
19 & 576 & 84 & 21 & 4 for each type\\
\bottomrule
\end{tabular}
\caption{Post rooms with sources $(0,0)$, $(2,0)$ and post $(1,0)$ (Theorem~\ref{thm:post}). ``Rooms'' counts the rooms with the sources and the post in these positions.}
\label{tab:post}
\end{table}

Type~3 of Figure~\ref{fig:post18} was found by the second author in 1996 but not published; the post room in his letter~\cite[Figure~3(c)]{Tok96} has 23 sides.

\subsection{Several posts}
A room on $\T$ \emph{with $k$ posts} is a union $R$ of tiles of $\T$ with connected interior whose boundary is a simple closed polygon, together with $k$ interior $B$- or $C$-points of $R$, the \emph{posts}, such that every other $B$- and $C$-point of $R$ is a corner. Light that reaches a post is absorbed, and the sources are the interior $A$-points. By Lemma~\ref{lem:blocking} the sources are invisible and pairwise mutually invisible. For $k=1$ these are the post rooms whose post is a $B$- or $C$-point; by the proof of Theorem~\ref{thm:post}, they include every post room with at most $19$ sides.

\begin{theorem}\label{thm:posts}
A room on $\T$ with at least two sources and $k$ posts has at least $17$ sides if $k=2$, and at least $16$ sides if $k=3$ or $k=4$. Up to congruence there are exactly six rooms with two sources, two posts and $17$ sides, and exactly four with two sources, three posts and $16$ sides \textup{(}Figures~\textup{\ref{fig:posts2}} and~\textup{\ref{fig:posts3}}\textup{)}. In each of them, after a symmetry of $\T$, the sources are $(0,0)$ and $(2,0)$, one post is $(1,0)$, and the other posts are $C$-points adjacent to a source.
\end{theorem}

\begin{proof}
Let $k_C$ and $k_B$ be the numbers of posts at $C$-points and at $B$-points. As in the proof of Lemma~\ref{lem:wings}(c), $T$ is the sum, over all $C$-points, of the number of tiles of $R$ there; this number is $0$, $1$ or $3$, and $4$ at a post. A $C$-point with three tiles has one full edge and a post at a $C$-point has two, so $T=3U+Q-2k_C$, where $Q$ is the number of $C$-points in one tile and $U=4M+f+\delta$ with $\delta\ge0$. The interior lattice points are the $M$ sources and the $k$ posts, and Lemma~\ref{lem:euler} gives $N+f=T-2(M+k)+2$. Hence
\begin{equation}\label{eq:posts}
N=10M+2+2f+3\delta+Q-4k_C-2k_B,
\end{equation}
which for $k=1$ is the equation in the proof of Theorem~\ref{thm:post}. For given bounds on $N$ and $k$, \eqref{eq:posts} gives $M=2$ and bounds $f$, and hence $T=N+2M+2k-2+f$: for $N\le17$ and $k\le2$, for $N\le16$ and $k\le3$, and for $N\le15$ and $k\le4$ it gives $T\le24$, $T\le27$ and $T\le29$. Take a shortest chain of tiles of $R$, each sharing an edge with the next, from the eight tiles at one source $s$ to the eight tiles at the other source $s'$. Its inner tiles lie outside the $16$ tiles at the two sources, so there are at most $T-16$ of them, and by Lemma~\ref{lem:chain}, with $r=1$, the sources $s$ and $s'$ have sup-distance at most $T-14$. This leaves finitely many tile sets. An exhaustive search over them, which allows $B$- and $C$-points in the interior and otherwise imposes the local conditions of a lattice room listed in the proof of Theorem~\ref{thm:tri}, finds no room with $k\le2$ and $N\le16$, none with $k\le4$ and $N\le15$, $24$ rooms with two posts and $17$ sides, and $9$ rooms with three posts and $16$ sides. With $k\le1$ and $N\le18$ it returns the $14$ rooms of Theorem~\ref{thm:post}. The types are counted as in the proof of Theorem~\ref{thm:post}; they are the orbits shown in Figures~\ref{fig:posts2} and~\ref{fig:posts3}.
\end{proof}

With at most four posts, the least number of sides is therefore $16$.

\subsection{One partition}
\begin{figure}[htbp]
\centering
\includegraphics[width=\textwidth]{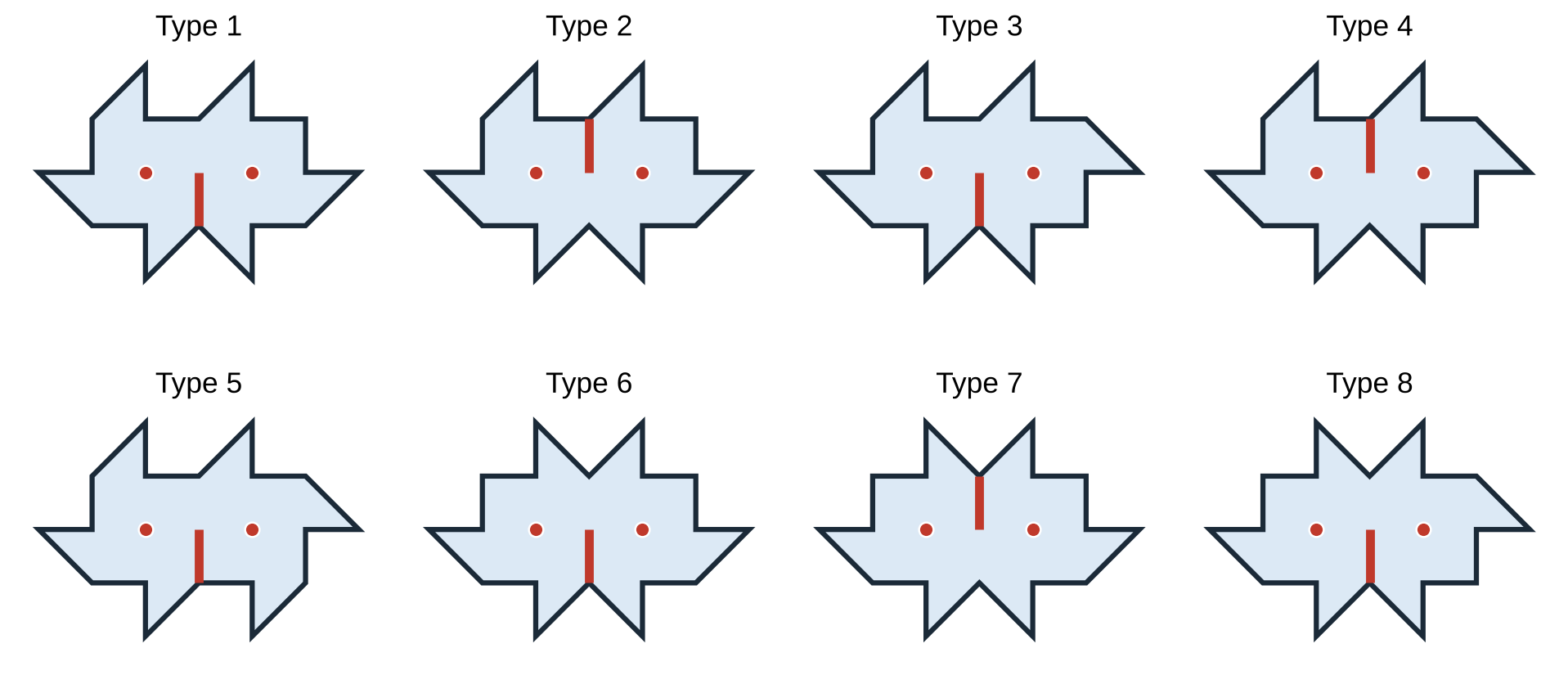}
\caption{The eight partition rooms with 20 sides (Corollary~\ref{cor:partition}). The partition is drawn as a thick segment. Type~5 is the room of~\cite[Figure~3(a)]{Tok96}.}
\label{fig:partition}
\end{figure}

\begin{corollary}\label{cor:partition}
A partition room with at least two sources has at least $20$ sides. Up to congruence there are exactly eight partition rooms with two sources and $20$ sides \textup{(}Figure~\textup{\ref{fig:partition}}\textup{)}. In each of them, after a symmetry of $\T$, the sources are $(0,0)$ and $(2,0)$, and the partition is the unit segment from $(1,0)$ to $(1,1)$ or to $(1,-1)$.
\end{corollary}

\begin{proof}
Let $R$ be a partition room with partition $J$ and free end $p$. Every $B$- and $C$-point of $R$ other than $p$ is a corner of $R$, so the relative interior of $J$ contains no such point. Removing $J$ and taking $p$ as a post therefore gives a post room with the same sources, whose number of sides is $N-2$. By Theorem~\ref{thm:post}, $N-2\ge18$. If $N=20$, the room without $J$ is one of the fourteen post rooms with $18$ sides and sources $(0,0)$, $(2,0)$, and $p=(1,0)$, since this point lies in four tiles. No diagonal mirror line passes through $(1,0)$, and a horizontal partition from $(1,0)$ would contain $(0,0)$ or $(2,0)$. So $J$ is vertical. In each of the fourteen rooms, the points $(1,1)$ and $(1,-1)$ are $B$-points of $R$, hence corners, and the open segments joining them to $(1,0)$ lie in the interior of $R$. The partition therefore runs from $(1,0)$ to $(1,1)$ or to $(1,-1)$, and both choices are allowed. This gives $28$ rooms. As in the proof of Theorem~\ref{thm:post}, the types are the orbits of the group of order four found there; there are eight, six of four rooms each and two of two rooms each.
\end{proof}

The partition room in the second author's letter~\cite[Figure~3(a)]{Tok96} is Type~5 of Figure~\ref{fig:partition}.

\subsection{One pinched hole}
\begin{figure}[htbp]
\centering
\includegraphics[width=\textwidth]{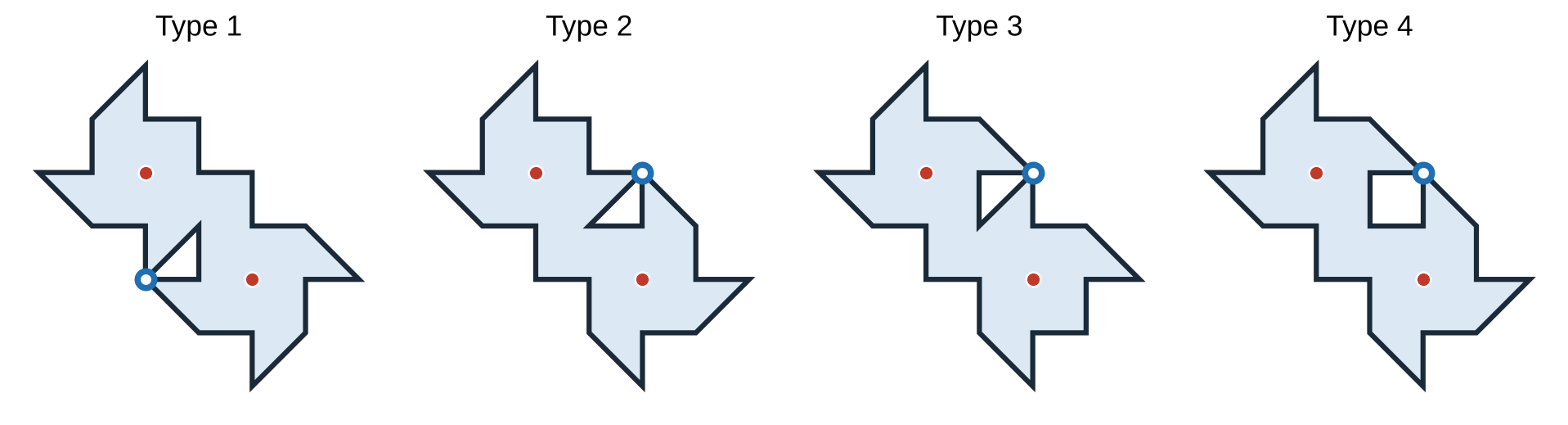}
\caption{The four pinched rooms with 22 sides (Theorem~\ref{thm:pinched}). The open circle marks the point $q$ where the hole meets the outer boundary. Type~2 is the room of~\cite[Figure~3(b)]{Tok96}, of Thomson~\cite{Thomson} and of~\cite[Figure~1]{Wolecki}, and Type~4 is Thomson's room with a square hole.}
\label{fig:pinched}
\end{figure}

\begin{theorem}\label{thm:pinched}
A pinched room with at least two sources has at least $22$ sides. Up to congruence there are exactly four pinched rooms with two sources and $22$ sides \textup{(}Figure~\textup{\ref{fig:pinched}}\textup{)}: three with a triangular hole and one with a square hole. In each of them, after a symmetry of $\T$, the sources are $(0,0)$ and $(2,-2)$.
\end{theorem}

\begin{proof}
Let $R$ be a pinched room with $M\ge2$ sources and $N\le22$ sides. Let $n_o$ and $n_h$ be the numbers of lattice points on the outer polygon and on the hole polygon, each counting $q$. Every lattice point on these polygons other than $q$ is a corner or a flat $A$-point, and $q$ counts as a corner of both, so $n_o+n_h=N+f$. The vertices and edges of the tiles of $R$ form a connected plane graph with $M+n_o+n_h-1$ vertices, $(3T+n_o+n_h)/2$ edges and $T+2$ faces: the tiles, the hole and the outer face. Euler's formula $(M+n_o+n_h-1)-(3T+n_o+n_h)/2+(T+2)=2$ gives $n_o+n_h=T-2M+2$, that is,
\begin{equation}\label{eq:pincheuler}
N=T-2M+2-f.
\end{equation}
The proof of Lemma~\ref{lem:wings}(c) applies at every $C$-point other than $q$. If $q$ is a $C$-point, then $R$ meets $q$ in two wedges, so $q$ lies in two opposite tiles of $R$ and in no full edge. Hence $T\ge12M+3f+2\varepsilon$, where $\varepsilon=1$ if $q$ is a $C$-point and $\varepsilon=0$ otherwise, and with~\eqref{eq:pincheuler} we get $N\ge10M+2+2f+2\varepsilon$. Thus $M=2$, $f=0$, $\varepsilon=0$, $N=22$ and $T=24$.

Each source lies in eight tiles and has four wings, as in Lemma~\ref{lem:wings}(b). The two sources have no common adjacent $C$-point, which would be interior, so these $24$ tiles are distinct and form all of $R$. The module of a source $s$ lies in the rotated square $|x-s_x|+|y-s_y|\le2$. Since the interior of $R$ is connected, the two modules share an edge, so these two squares share a segment of positive length. After a symmetry of $\T$ the sources are $(0,0)$ and $(2a,2b)$ with $0\le b\le a$, and the squares share such a segment only for $(a,b)=(1,0)$ and $(1,1)$. For $(1,0)$ the sources share the adjacent $C$-point $(1,0)$, which is impossible. So, after a reflection, the sources are $(0,0)$ and $(2,-2)$.

The eight wings can be chosen independently, two ways each, which gives $256$ candidates. We tested each candidate, with exact integer arithmetic, for connected interior, for a boundary consisting of two simple polygons that meet in exactly one point, for the corner condition away from that point and for the number of sides. Fourteen rooms pass these tests. The isometries that map $\{(0,0),(2,-2)\}$ onto itself are the identity, the reflections in the lines $x+y=0$ and $x-y=2$ and the rotation through $180\dg$ about $(1,-1)$; all four are symmetries of $\T$, and by Lemma~\ref{lem:congruence} every congruence between two of the rooms is one of them. The fourteen rooms form four types: three of four rooms each, with triangular holes, and one of two rooms, with a square hole.
\end{proof}

\begin{remark}\label{rem:pinchcount}
In the three rooms with a triangular hole the outer polygon turns at $q$; it has $19$ corners and the hole has $3$. In the room with a square hole, Type~4 of Figure~\ref{fig:pinched}, the outer polygon is straight at $q$. Near $q$ the hole occupies a right angle inside the straight angle of the outer boundary, and the room consists of two wedges of $45\dg$, each bounded by a piece of the outer boundary and a side of the hole. Light reaching $q$ is absorbed, and the two pieces of the outer boundary bound different wedges, so we count them as two sides. This gives $18+4=22$. If one counts only maximal straight segments of the boundary, this room has $17+4=21$ sides. Four sides end at $q$, and a line through $q$ contains at most two of them, so this count is at least $N-2$. In the four types of Figure~\ref{fig:pinched} no side of the hole is collinear with a side of the outer polygon at $q$, so their counts are $22$, $22$, $22$ and $21$. By Theorem~\ref{thm:pinched}, a pinched room with two sources therefore has at least $21$ maximal straight boundary segments, and among the rooms with $22$ sides Type~4 is the only one with $21$.
\end{remark}

\begin{remark}\label{rem:wolecki}
Type~2 of Figure~\ref{fig:pinched} is the room with a triangular hole in the second author's letter~\cite[Figure~3(b)]{Tok96}, and it is also the 22-sided room of Wolecki~\cite[Figure~1]{Wolecki}; Thomson~\cite{Thomson} found it as well. Type~4 is the room with a square hole that Thomson~\cite{Thomson} sent to the second author in February 2001. Thomson counted its sides as maximal straight segments, found $21$, and noted that the room is symmetric in the line through the diagonal of the hole, and that folding along this line gives ``the smallest self-unilluminable room''. Each half is a pinwheel, so Theorem~\ref{thm:pinwheel} confirms this for hole-free lattice rooms on $\T$. We matched the coordinates of these rooms with our list by computer.
\end{remark}

\subsection{One hole}
A hole that does not touch the outer boundary costs more sides.

\begin{theorem}\label{thm:hole}
A lattice room on $\T$ with exactly one hole and at least two sources has at least $26$ sides. Up to congruence there are exactly $12$ such rooms with $26$ sides \textup{(}Figure~\textup{\ref{fig:hole}}\textup{)}. In each of them, after a symmetry of $\T$, the sources are $(0,0)$ and $(2,2)$.
\end{theorem}

\begin{proof}
With $h=1$, \eqref{eq:euler} gives $N=T-2M-f$, and $T=3U+Q$ with $U=4M+f+\delta$ as in the proof of Lemma~\ref{lem:wings}(c). Hence $N=10M+2f+3\delta+Q$. Let $N\le26$. Then $M=2$, $N=20+2f+3\delta+Q$ and $T=N+4+f$, so $f\le3$ and $T\le33$. The two modules are $24$ distinct tiles, since two sources have no common adjacent $C$-point; so at most nine tiles lie outside them. A shortest chain of tiles of $R$, each sharing an edge with the next, from one module to the other has its $j\le9$ inner tiles outside the modules, so by Lemma~\ref{lem:chain}, with $r=2$, the sources have sup-distance at most $j+4\le13$. After a symmetry of $\T$ the sources are $(0,0)$ and $(2a,2b)$ with $0\le b\le a\le6$. An exhaustive search over at most $33$ tiles in each of these $27$ cases, with the local conditions of a lattice room and a boundary of two components, finds no room with fewer than $26$ sides and $42$ rooms with $26$ sides. All of them have sources $(0,0)$ and $(2,2)$, and they form $12$ orbits of the symmetries of $\T$ that map the pair of sources onto itself. By Lemma~\ref{lem:congruence}, these orbits are the types.
\end{proof}

The rooms of Theorems~\ref{thm:pinched} and~\ref{thm:hole} have the same pair of sources up to a symmetry of $\T$. Moving the hole away from the outer boundary costs four sides.

\section{Hole-free rooms with two sources}\label{sec:holefree}

\begin{theorem}\label{thm:tri}
A hole-free lattice room on $\T$ with at least two sources has at least $24$ sides. Up to congruence there are exactly $12$ such rooms with $24$ sides and exactly $89$ with $25$ sides \textup{(}Figures~\textup{\ref{fig:t24}--\ref{fig:t25o}}\textup{)}. Each has two sources. After a symmetry of $\T$ the sources are $(0,0)$ and $(4,0)$, except in $24$ of the types with $25$ sides, where they are $(0,0)$ and $(4,2)$.
\end{theorem}

\begin{proof}
Let $R$ be such a room with $N\le25$. By Corollary~\ref{cor:bound}, $N\ge10M+2+2f$, so $M=2$ and $f\le1$, and~\eqref{eq:euler} gives $T=N+2+f\le28$. By Lemma~\ref{lem:wings}(b), $R$ contains the eight tiles at each source and the four wings of each source. Two sources have no common adjacent $C$-point, since such a point would lie in four tiles, and it follows that these $24$ tiles are distinct.

After a symmetry of $\T$ the sources are $(0,0)$ and $(2a,2b)$ with $0\le b\le a$. Every side of $R$ has horizontal extent at most $2$: its endpoints are lattice points, its relative interior contains no $B$- or $C$-point, and on a mirror line that is not vertical the lattice points follow each other at horizontal distance $1$, with at most every second one an $A$-point. The boundary of $R$ is a closed curve, so the total horizontal extent of its sides, at most $2N$, is at least twice the width of $R$. The tiles at the sources span the strip $-1\le x\le2a+1$, so $2a+2\le N\le25$ and $a\le11$. This leaves $\sum_{a=1}^{11}(a+1)=77$ positions for the second source.

For each position, a computer search listed every set of at most $28$ tiles that contains the sixteen tiles at the two sources, has connected interior and meets the local conditions of a lattice room at every lattice point: at a $C$-point it has $0$, $1$ or $3$ tiles; at an $A$-point its tiles are absent, are all eight tiles there, or form one block of consecutive tiles; and at a $B$-point they are absent or form one block of consecutive tiles that makes a corner. The search starts from the sixteen tiles at the sources and adds tiles to them. At a lattice point where the current tiles violate the local conditions, it branches over all minimal sets of tiles at that point whose addition satisfies them; every room containing the current set contains one of these. Otherwise it adds, in every possible way, one tile that shares an edge with the current set; every larger set with connected interior that contains the current set contains such a tile. Repeated sets are discarded. The search region consists of the tiles in the square $[-30,31]^2$. A candidate has at most $28$ tiles, contains the tiles at the origin and has connected interior, so each of its tiles is joined to those tiles by a chain of fewer than $28$ tiles, and the argument of Lemma~\ref{lem:chain} shows that it lies in $[-28,28]^2$; so no candidate is cut off. Finally we kept the sets whose boundary is one simple polygon satisfying the corner condition, with at most $25$ sides. Only the positions $(4,0)$ and $(4,2)$ give rooms, and the results are in Table~\ref{tab:tri}.

By Lemma~\ref{lem:congruence}, a congruence between two of these rooms maps $\Z^2$ onto itself and sources onto sources. The isometries with these two properties that map $\{(0,0),(4,0)\}$ onto itself are the identity, the reflections in $y=0$ and $x=2$ and the rotation through $180\dg$ about $(2,0)$; for $\{(0,0),(4,2)\}$ they are the identity and the rotation through $180\dg$ about $(2,1)$. All of them are symmetries of $\T$, and the types in Table~\ref{tab:tri} are the orbits of these two groups.
\end{proof}

\begin{table}[htbp]
\centering\small
\begin{tabular}{@{}cccccl@{}}
\toprule
Sides & Second source & Rooms & Types & Tiles & Labels in Appendix~\ref{app:plates}\\
\midrule
24 & $(4,0)$ & 48 & 12 & 27 & T24-001 to T24-012\\
25 & $(4,0)$ & 260 & 65 & 28 & T25-001 to T25-065\\
25 & $(4,2)$ & 48 & 24 & 27 or 28 & T25-066 to T25-089\\
\bottomrule
\end{tabular}
\caption{Hole-free lattice rooms on $\T$ with two sources and at most 25 sides (Theorem~\ref{thm:tri}). The first source is $(0,0)$; ``Rooms'' counts the rooms with the sources in these positions. For the pair $(0,0)$, $(4,2)$, half of the rooms have 27 tiles and half have 28.}
\label{tab:tri}
\end{table}

Castro~\cite{Castro}, the second author~\cite[Figure~2]{Tok96} and Thomson~\cite{Thomson} found rooms with 24 sides; they are T24-003, T24-008 and T24-012. The second author's room is room~b of~\cite{Stewart97}, where it is conjectured to be the smallest unilluminable room, and Theorem~\ref{thm:tri} shows that no hole-free lattice room on $\T$ with two sources has fewer sides. The straight snake with two sources of Section~\ref{sec:many} is Castro's room T24-003. The rooms with 25 sides are the smallest rooms in this class with an odd number of sides. The first unilluminable room with an odd number of sides had 27 sides; it is in the second author's letter~\cite[Figure~5]{Tok96}, see also room~d of~\cite{Stewart97}. Wolecki's 22-sided room has a pinched hole (Section~\ref{sec:devices}), and so does not belong to this class.

\begin{theorem}\label{thm:sq}
A hole-free lattice room on $\Sq$ with at least two sources has at least $34$ sides. Up to congruence there are exactly $14$ such rooms with $34$ sides \textup{(}Figure~\textup{\ref{fig:s34}}\textup{)}. Each has two sources, $19$ squares and two flat points.
\end{theorem}

\begin{proof}
Let $R$ be such a room with $N\le34$. As in the proof of Theorem~\ref{thm:square}, $4T\le3N-6M-6$, and $M\ge2$ gives $T\le21$. The boundary of $R$ has $N/2$ horizontal sides, each of length at most $2$ by the corner condition. Their total length is at least twice the width of $R$, as in the proof of Theorem~\ref{thm:tri}, so the width of $R$ is at most $N/2$, which is at most $17$. With two of the sources at $(0,0)$ and $(2a,2b)$, $0\le b\le a$, the squares at these two sources span the strip $-1\le x\le2a+1$, so $a\le7$, which leaves $\sum_{a=1}^{7}(a+1)=35$ positions. The search described in the proof of Theorem~\ref{thm:tri}, run with up to $22$ squares, gives Table~\ref{tab:sq}. A second program, which grows connected sets of squares one square at a time, finds the same $14$ types with up to $21$ squares. Each of the $40$ rooms has $19$ squares, two sources and two flat points, so $N=2\cdot19-4+2-2=34$ by~\eqref{eq:euler}. By Lemma~\ref{lem:congruence}, the types are the orbits of the symmetry groups listed in Table~\ref{tab:sq}.
\end{proof}

Tokarsky's rectilinear room~\cite[Figure~9]{Tok95} has 34 sides and is S34-001, so Theorem~\ref{thm:sq} shows that it is smallest in its class.

\begin{table}[htbp]
\centering\small
\begin{tabular}{@{}cccccl@{}}
\toprule
Second source & Rooms & Symmetry group & Types & Squares & Labels in Appendix~\ref{app:plates}\\
\midrule
$(4,0)$ & 24 & order 4 & 6 & 19 & S34-001 to S34-006\\
$(4,2)$ & 14 & order 2 & 7 & 19 & S34-007 to S34-013\\
$(4,4)$ & 2 & order 4 & 1 & 19 & S34-014\\
\midrule
Total & 40 & & 14 & &\\
\bottomrule
\end{tabular}
\caption{Hole-free lattice rooms on $\Sq$ with two sources and 34 sides (Theorem~\ref{thm:sq}). The first source is $(0,0)$; the symmetry group consists of the symmetries of $\Sq$ that map the pair of sources onto itself.}
\label{tab:sq}
\end{table}

\section{Rooms with many sources}\label{sec:many}

Lattice rooms with several sources can be built from modules, each consisting of the tiles at a source and its four wings, joined by small pieces.

\begin{definition}\label{def:assembly}
Let $s$ be an $A$-point. A \emph{module} at $s$ is the union of the eight tiles of $\T$ at $s$ and, at each of the four $C$-points adjacent to $s$, one of the two tiles there that do not contain $s$. There are sixteen modules at $s$, and each is a disc of twelve tiles. A \emph{connector} is the union of three consecutive tiles of $\T$ at a $C$-point; it is a disc of three tiles. An \emph{assembly} is a room $R$ that is the union of $n\ge1$ modules at distinct points $s_1,\dots,s_n$ and $c\ge0$ connectors, any two of which have disjoint interiors, such that:
\begin{enumerate}[label=\textup{(\roman*)}]
\item each connector $K$ meets exactly two modules, in one edge each; these two edges have a common endpoint $q_K$, which is an $A$-point; the two modules meet each other exactly in $q_K$; and no other two of the modules and connectors meet;
\item the graph $G$ whose vertices are the modules, with one edge for each connector joining the two modules it meets, is connected;
\item $R$ is a lattice room, its sources are $s_1,\dots,s_n$, and its flat points are the points $q_K$.
\end{enumerate}
\end{definition}

\begin{theorem}\label{thm:assembly}
An assembly with $n$ modules, $c$ connectors and $h$ holes satisfies $c=n-1+h$ and has exactly $12n$ sides. Its sources are invisible and pairwise mutually invisible. In particular, if $G$ is a tree then $R$ has no holes and $n-1$ connectors, and if $G$ is a cycle then $R$ has one hole and $n$ connectors.
\end{theorem}

\begin{proof}
The modules and connectors have disjoint interiors, so $T=12n+3c$, and by condition~(iii) $M=n$ and $f=c$. For a set $X$ that is a union of tiles, tile edges and tile vertices of $\T$, let $\chi(X)=V-E+F$, where $V$, $E$ and $F$ are the numbers of tile vertices, tile edges and tiles contained in $X$. Each of these three numbers obeys inclusion and exclusion over unions of such sets, and so does $\chi$. A disc made of tiles has $\chi=1$, by Euler's formula as in the proof of Lemma~\ref{lem:euler}, and so does a single edge or a single point. The same proof gives $\chi(R)=1-h$. The modules and connectors are discs. By condition~(i), the nonempty intersections of two of them are, for each connector $K$, the two edges where $K$ meets its modules and the point $q_K$ where these modules meet. The only nonempty intersections of three of them are the points $q_K$, one for each connector. Inclusion and exclusion give $\chi(R)=(n+c)-3c+c=n-c$, and with $\chi(R)=1-h$ we get $h=c-n+1$. Now~\eqref{eq:euler} gives
\[
N=(12n+3c)-2n+2-2(c-n+1)-c=12n.
\]
The graph $G$ is connected with $n$ vertices and $c$ edges, so its cycle rank $c-n+1$, the number of its edges outside a spanning tree, equals $h$, which gives the statements about trees and cycles. The sources are invisible and mutually invisible by Corollary~\ref{cor:blocking}.
\end{proof}

\begin{figure}[htbp]
\centering
\includegraphics[width=\textwidth]{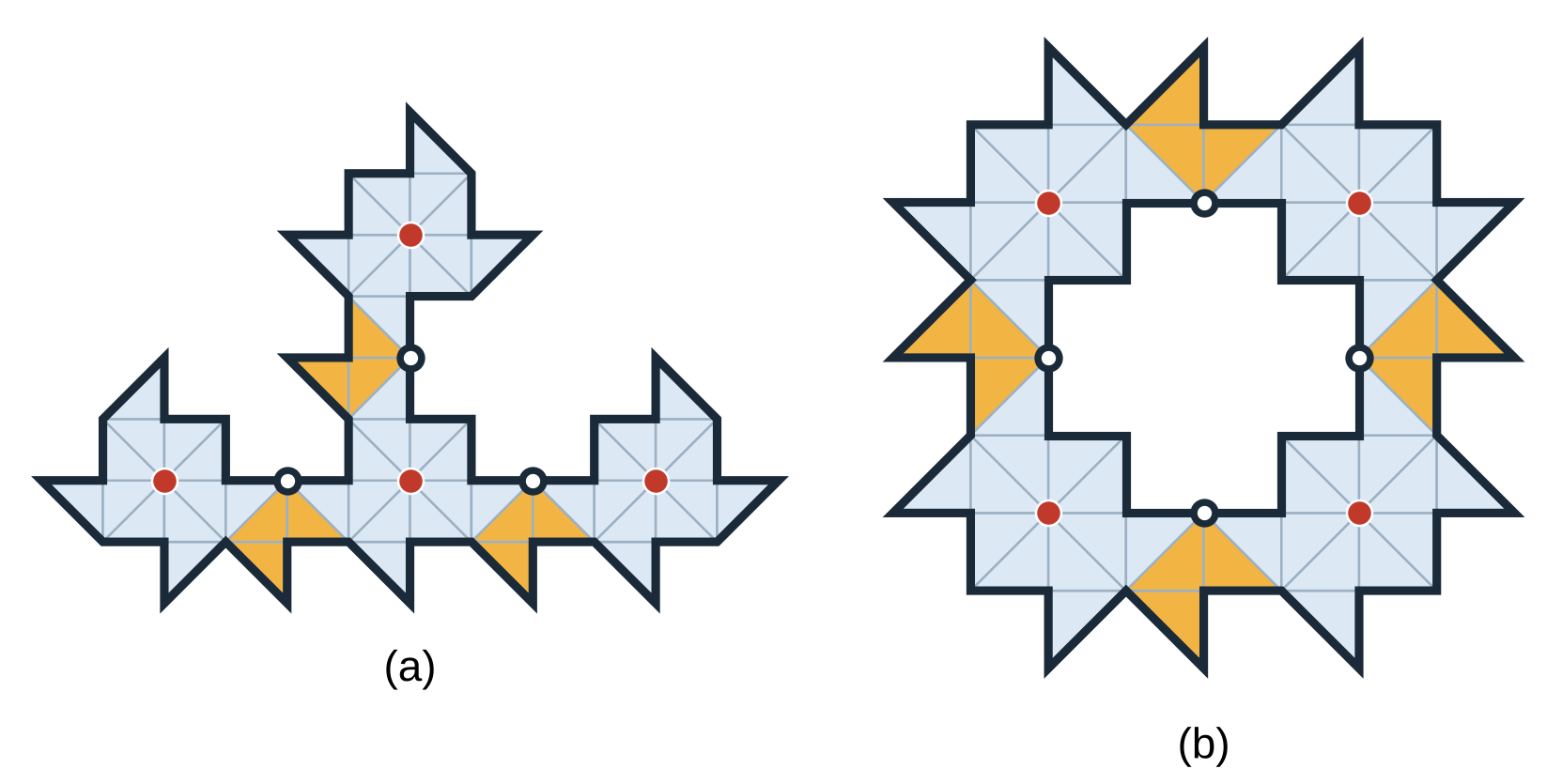}
\caption{Two assemblies with four sources (dots) and 48 sides. Modules are light blue, connectors are orange, and the open circles are the flat points where connectors are attached. (a)~A tree with three connectors. (b)~A cycle with four connectors and one hole.}
\label{fig:many}
\end{figure}

\begin{example}\label{ex:tree}
Figure~\ref{fig:many}(a) shows an assembly whose graph is a tree with a vertex of degree three. Its sources are $(-4,0)$, $(0,0)$, $(4,0)$ and $(0,4)$; it has three connectors, $57$ tiles, three flat points and $48$ sides, and its corners are $13$ $A$-points, $16$ $B$-points and $19$ $C$-points. Figure~\ref{fig:many}(b) shows an assembly whose graph is a cycle of length four, with sources $(0,0)$, $(4,0)$, $(0,4)$, $(4,4)$. It has four connectors, $60$ tiles, four flat points and one hole; its outer boundary has $36$ sides and the hole $12$, and its corners are $12$ $A$-points, $16$ $B$-points and $20$ $C$-points. Compared with the tree, the cycle has one more connector and one hole, and the same number of sides, as Theorem~\ref{thm:assembly} predicts.
\end{example}

The simplest assemblies are chains of modules along a line. Appendix~\ref{app:snakes} gives explicit vertex lists for two infinite families: the \emph{straight snakes} $\Gamma_n$, $n\ge2$, with sources $(4j,0)$ for $0\le j<n$, and the \emph{bent snakes} $\Gamma'_n$, $n\ge3$, with sources $(4j,0)$ for $0\le j\le n-2$ together with $(0,4)$.

\begin{proposition}\label{prop:snakes}
For every $n\ge2$ the straight snake $\Gamma_n$, and for every $n\ge3$ the bent snake $\Gamma'_n$, is a hole-free lattice room on $\T$ with $n$ sources, $n-1$ flat points, $15n-3$ tiles and $12n$ sides. Its corners are $3n+1$ $A$-points, $4n$ $B$-points and $5n-1$ $C$-points.
\end{proposition}

The proof, for all $n$, is in Appendix~\ref{app:snakes}. The tile count $15n-3=12n+3(n-1)$ is that of an assembly of $n$ modules and $n-1$ connectors along a path, and $\Gamma_2$ is the 24-sided room T24-003.

\begin{theorem}\label{thm:three}
A hole-free lattice room on $\T$ with at least three sources has at least $36$ sides.
\end{theorem}

\begin{proof}
If $M\ge4$, Corollary~\ref{cor:bound} gives $N\ge42$. Let $M=3$ and $N\le35$. As in the proof of Theorem~\ref{thm:hole}, now with $h=0$, $N=32+2f+3\delta+Q$ and $T=36+3f+3\delta+Q$, so $2f+3\delta+Q\le3$ and $T\le40$. The three modules are $36$ distinct tiles, so at most four tiles of $R$ lie outside them. Since the interior of $R$ is connected, a shortest chain of tiles of $R$ from the module of one source to the union of the other two modules ends at the module of a second source, and a shortest chain from the module of the third source to the union of the first two modules ends at one of them. The inner tiles of each chain lie outside the modules, so there are at most four of them. By Lemma~\ref{lem:chain}, with $r=2$, each chain joins two sources at sup-distance at most $8$. Up to symmetries of $\T$ this leaves $281$ positions of three sources, and an exhaustive search over at most $40$ tiles finds no room in any of them. The snakes $\Gamma_3$ and $\Gamma'_3$ have $36$ sides.
\end{proof}

For $n\le3$, the least number of sides of a hole-free lattice room on $\T$ with $n$ sources is therefore $12n$.

\section{Relaxing the corner condition}\label{sec:relax}

Lemma~\ref{lem:blocking} needs only the points of $\Lambda\setminus2\Lambda$ to be corners. On $\T$, with $\Lambda=\Z^2$, these are all $B$- and $C$-points. On other tilings $\Lambda\setminus2\Lambda$ can be a proper subset of the tile vertices, and then fewer corners are needed. A symmetry of the room can also take the place of a corner.

\subsection{Rooms with 27 sides on the tiling \texorpdfstring{$\K$}{K}}\label{sec:27}
Let $\K$ be the tiling of the plane by triangles with angles $30\dg$, $60\dg$, $90\dg$ generated by reflections in their sides, with short leg of length $1$. Its vertices are of three kinds: \emph{$A$-points}, where twelve tiles meet with their $30\dg$ angles; \emph{$B$-points}, where six tiles meet; and \emph{$C$-points}, where four tiles meet with their right angles. The $A$-points form a triangular lattice $\Lambda_A$ with minimal distance $2\sqrt3$, spanned by two vectors $e_1$ and $e_2$ of length $2\sqrt3$ at an angle of $60\dg$. Let $\Lambda=\tfrac12\Lambda_A$ be the lattice spanned by $e_1/2$ and $e_2/2$. Then $2\Lambda=\Lambda_A$, and $\Lambda\setminus2\Lambda$ consists of the points $p+e_1/2$, $p+e_2/2$ and $p+(e_1-e_2)/2$ with $p\in\Lambda_A$, which are the midpoints of pairs of nearest $A$-points; these are the $C$-points. A reflection in a mirror line of $\K$, a line that contains tile edges, maps the tiling onto itself, so it maps the $A$-points onto the $A$-points and the $C$-points onto the $C$-points, and therefore $\Lambda$ onto $\Lambda$ and $2\Lambda$ onto $2\Lambda$. By Lemma~\ref{lem:blocking}, if every $C$-point of a union $R$ of tiles of $\K$ is a corner of $R$, then the $A$-points of $R$ are invisible and pairwise mutually invisible. We again call the interior $A$-points the \emph{sources} of $R$. The $B$-points of $R$ may be interior points or flat points.

We place an $A$-point at the origin, with the $A$-point $(2\sqrt3,0)$ as one of its nearest neighbours, and use coordinates in which $(u,v)$ denotes the point $(\sqrt3\,u/2,\,v/2)$. Then $\Lambda_A$ is spanned by $(4,0)$ and $(2,6)$, and the $C$-points adjacent to an $A$-point are the six $C$-points at distance $\sqrt3$ from it.

\begin{proposition}\label{prop:27}
The $27$-gon of Figure~\textup{\ref{fig:relax}(a)}, whose vertices are listed in Appendix~\textup{\ref{app:coords}}, is a hole-free union of $42$ tiles of $\K$ in which every $C$-point is a corner. Its interior $A$-points $P=(0,0)$ and $Q=(8,0)$ are invisible and mutually invisible. It has $8$ corners of angle $30\dg$, $2$ of angle $60\dg$, $3$ of angle $120\dg$ and $14$ of angle $270\dg$, and the corners of angle $270\dg$ are exactly its $C$-points. Two of its $B$-points, $(2,-2)$ and $(6,-2)$, are interior points, and eight lie in the relative interiors of sides.
\end{proposition}

\begin{proof}
The geometric statements were verified with exact arithmetic, and the invisibility follows from Lemma~\ref{lem:blocking}.
\end{proof}

\begin{theorem}\label{thm:min27}
Let $R$ be a union of tiles of $\K$ whose boundary is a simple closed polygon, with at least two interior $A$-points, in which every $C$-point is a corner. Then $R$ has at least $27$ sides. Up to congruence there are exactly $15$ such rooms with $27$ sides \textup{(}Figure~\textup{\ref{fig:k27}}\textup{)}; each consists of $42$ tiles, and after a symmetry of $\K$ its sources are $P=(0,0)$ and $Q=(8,0)$.
\end{theorem}

\begin{proof}
Suppose that $R$ has $N\le26$ sides, and let $r$ be the number of its corners of angle $270\dg$. Every $C$-point of $R$ is a corner with angle $90\dg$ or $270\dg$, and every corner has angle at least $30\dg$. The exterior angles of a simple polygon add up to $360\dg$, so
\[
360\le150(N-r)-90r,\qquad\text{that is,}\qquad N\ge\frac{12+8r}{5},
\]
and $r\le14$. At an interior $A$-point $s$, all twelve tiles at $s$ lie in $R$, and they supply two of the four tiles at each of the six $C$-points adjacent to $s$. Since these $C$-points are corners, each lies in a third tile and has angle $270\dg$. Two interior $A$-points have no common adjacent $C$-point, since it would be interior. Hence $R$ has exactly two interior $A$-points $P$ and $Q$, which account for $12$ of the at most $14$ corners of angle $270\dg$.

For a $C$-point $c$ of $R$, call the set of tiles of $R$ at $c$ the \emph{cell} of $c$. Every tile has exactly one $C$-vertex, its right angle, so the cells partition the tiles of $R$, and each cell has one or three tiles because its $C$-point is a corner. Tiles in different cells can share only an edge joining an $A$-point to a $B$-point, the hypotenuse, and the $C$-points of two tiles with a common hypotenuse are mirror images in it, at distance $\sqrt3$. A cell with one tile has only one hypotenuse. In the graph whose vertices are the cells, with an edge between two cells that share an edge of a tile, every cell with one tile is therefore a leaf. This graph is connected, because the interior of $R$ is connected, and it stays connected when the leaves are removed. The cells with three tiles thus form a connected graph, whose edges join $C$-points at distance $\sqrt3$. A shortest path in it from the six cells around $P$ to the six cells around $Q$ passes through at most two other cells, so it has at most three edges, and $|PQ|\le\sqrt3+3\sqrt3+\sqrt3=5\sqrt3$. After a symmetry of $\K$, $Q-P=m(4,0)+n(2,6)$ with $m\ge n\ge0$, and $|PQ|^2=12(m^2+mn+n^2)\le75$ leaves $(m,n)=(1,0)$, $(1,1)$ and $(2,0)$. For $(1,0)$ the two points share an adjacent $C$-point. For $(1,1)$ and $(2,0)$ we list the possible sets of $C$-points with three tiles: sets of at most $14$ $C$-points that contain the twelve $C$-points adjacent to $P$ and $Q$ and are connected when $C$-points at distance $\sqrt3$ are joined. There are $13$ such sets: $11$ for $(1,1)$ and $2$ for $(2,0)$. In each case every other cell of $R$ has one tile and is adjacent to one of these cells, so only finitely many sets of tiles are possible. An exhaustive search over them requires at every vertex that the tiles of $R$ there are absent, are all present (only at an $A$- or $B$-point), or form one block of consecutive tiles; it requires three tiles at the $C$-points of the set and at most one at every other $C$-point, and exterior angles that add up to $360\dg$. It finds no room with at most $26$ corners in any of the $13$ cases. For $N=27$ the same inequality gives $r\le15$, so at most three cells with three tiles are not forced, $|PQ|\le6\sqrt3$, and $(m,n)$ is one of $(1,1)$, $(2,0)$, $(2,1)$ and $(3,0)$. There are $78$, $24$, $1$ and $0$ such sets of at most $15$ $C$-points, and the search over them lists $50$ rooms with $27$ corners. All have $(m,n)=(2,0)$, and they form $15$ types. We separated the types with a complete congruence invariant of a polygon: the matrix of squared distances between its corners, taken in cyclic order, up to rotation and reversal of that order.
\end{proof}

By Proposition~\ref{prop:27}, the bound of Theorem~\ref{thm:min27} is attained.

\subsection{A nonagon with three invisible corners}
\begin{proposition}\label{prop:nonagon}
In coordinates in which $(u,v)$ denotes the point $(u+v/2,\,\sqrt3\,v/2)$, let $R$ be the nonagon with vertices
\[
(0,0),\ (4,-2),\ (3,0),\ (6,0),\ (4,4),\ (3,3),\ (0,6),\ (-2,4),\ (0,3)
\]
\textup{(}Figure~\textup{\ref{fig:relax}(b))}. Its corners $A_1=(0,0)$, $A_2=(6,0)$ and $A_3=(0,6)$ are invisible and pairwise mutually invisible.
\end{proposition}

The nonagon, found by the second author, is the equilateral triangle $A_1A_2A_3$ with a right triangle with angles $30\dg$ and $60\dg$ attached along half of each side. It is a union of nine tiles of a tiling by triangles with angles $30\dg$, $60\dg$ and $90\dg$, with $A$-points at $A_1$, $A_2$, $A_3$ and $C$-points at the midpoints $(3,0)$, $(3,3)$, $(0,3)$ of the sides of the triangle. Its corners have angles $90\dg$ at the $A_i$, $60\dg$ at the tips of the attached triangles and $270\dg$ at the midpoints.

\begin{proof}
Take for $\Lambda$ the lattice $3\Z^2$ in these coordinates, so that $2\Lambda=6\Z^2$. In order from the side through $(0,0)$ and $(4,-2)$, the reflections in the lines of the nine sides are
\begin{gather*}
(u,v)\mapsto(u,\,-u-v),\quad(6-u-v,\,v),\quad(u+v,\,-v),\quad(12-u-v,\,v),\quad(v,\,u),\\
(6-v,\,6-u),\quad(v-6,\,u+6),\quad(u,\,6-u-v),\quad(-u,\,u+v).
\end{gather*}
Each has a linear part with integer entries and a translation part in $6\Z^2$, so it maps $\Lambda$ onto $\Lambda$ and $2\Lambda$ onto $2\Lambda$, and condition~(i) of Lemma~\ref{lem:blocking} is satisfied. The points of $\Lambda\setminus2\Lambda$ in $R$ are $(3,0)$, $(3,3)$ and $(0,3)$, which are corners of $R$, so condition~(ii) is satisfied. The points of $2\Lambda$ in $R$ are $A_1$, $A_2$ and $A_3$, and Lemma~\ref{lem:blocking} applies.
\end{proof}

The nonagon has no interior $A$-point, and the centre $(2,2)$ of the triangle, a $B$-point of the tiling, is an interior point.

\subsection{A symmetry in place of a corner}

\begin{theorem}\label{thm:symmetry}
Let $\Lambda$ be a lattice and let $R$ be a room that satisfies condition~\textup{(i)} of Lemma~\textup{\ref{lem:blocking}}. Let $\sigma$ be a reflection or a rotation through $180\dg$ with $\sigma(R)=R$, and let $U$ and $V$ be disjoint nonempty finite sets of interior points of $R$ in $2\Lambda$ with $\sigma(U)=U$ and $\sigma(V)=V$. Suppose that every point of $\Lambda\setminus2\Lambda$ that lies in $R$ is a corner of $R$, except for finitely many points $F$ fixed by $\sigma$, each of which is either an interior point of $R$ such that $\sigma$ is the rotation through $180\dg$ about $F$, or a point in the relative interior of a side such that $\sigma$ is the reflection in the line through $F$ perpendicular to that side. Then no light ray joins a point of $U$ to a point of $V$.
\end{theorem}

\begin{proof}
Suppose that such rays exist. Unfolding shows that the squared length of a ray from a point of $2\Lambda$ to a point of $2\Lambda$ is the squared length of a vector of $2\Lambda$; these values form a discrete set. Choose a shortest ray $\alpha$ from a point of $U$ to a point of $V$, parametrised by arc length on $[0,\ell]$. The proof of Lemma~\ref{lem:blocking} gives a time $t\in(0,\ell)$ at which $\alpha$ passes through a point of $\Lambda\setminus2\Lambda$. This point is not a corner, so it is one of the points $F$. If $F$ is interior, $\alpha$ passes straight through $F$, and $\sigma$, the rotation through $180\dg$ about $F$, maps the incoming segment onto the outgoing one with its direction reversed. If $F$ lies on a side, $\alpha$ reflects at $F$, and $\sigma$, the reflection in the normal line at $F$, again maps the incoming segment onto the outgoing one with its direction reversed. Since $\sigma$ maps light rays to light rays and a light ray is determined by any of its segments,
\[
\alpha(t+u)=\sigma\bigl(\alpha(t-u)\bigr)\qquad\text{for }0\le u\le\min(t,\ell-t).
\]
If $2t>\ell$, then $\alpha(2t-\ell)=\sigma(\alpha(\ell))\in V$, and $\alpha$ restricted to $[0,2t-\ell]$ is a shorter ray from $U$ to $V$. If $2t=\ell$, then $\alpha(\ell)=\sigma(\alpha(0))\in U$, which contradicts $U\cap V=\emptyset$. If $2t<\ell$, then $\alpha(2t)=\sigma(\alpha(0))\in U$, and $\alpha$ restricted to $[2t,\ell]$ is a shorter ray from $U$ to $V$. Each case contradicts the choice of $\alpha$.
\end{proof}

\begin{example}\label{ex:42}
Let $R$ be the 42-gon of Figure~\ref{fig:relax}(c), a union of tiles of $\T$ whose vertices are listed in Appendix~\ref{app:coords}. The rotation $\sigma(x,y)=(-x,-2-y)$ through $180\dg$ about $F=(0,-1)$ maps $R$ onto itself. With $\Lambda=\Z^2$, the only $B$- or $C$-point of $R$ that is not a corner is $F$, an interior $C$-point. The sets $U=\{P,P'\}=\{(0,0),(0,-2)\}$ and $V=\{Q,Q'\}=\{(4,0),(-4,-2)\}$ consist of interior $A$-points and are mapped onto themselves by $\sigma$. By Theorem~\ref{thm:symmetry}, no light ray joins a point of $U$ to a point of $V$; in particular $P$ and $Q$ are mutually invisible. The theorem says nothing about pairs within $U$, and indeed the segment from $P$ to $P'$ passes through $F$ and is a light ray; it also says nothing about whether the points of $U$ and $V$ are invisible.
\end{example}

\begin{figure}[htbp]
\centering
\includegraphics[width=\textwidth]{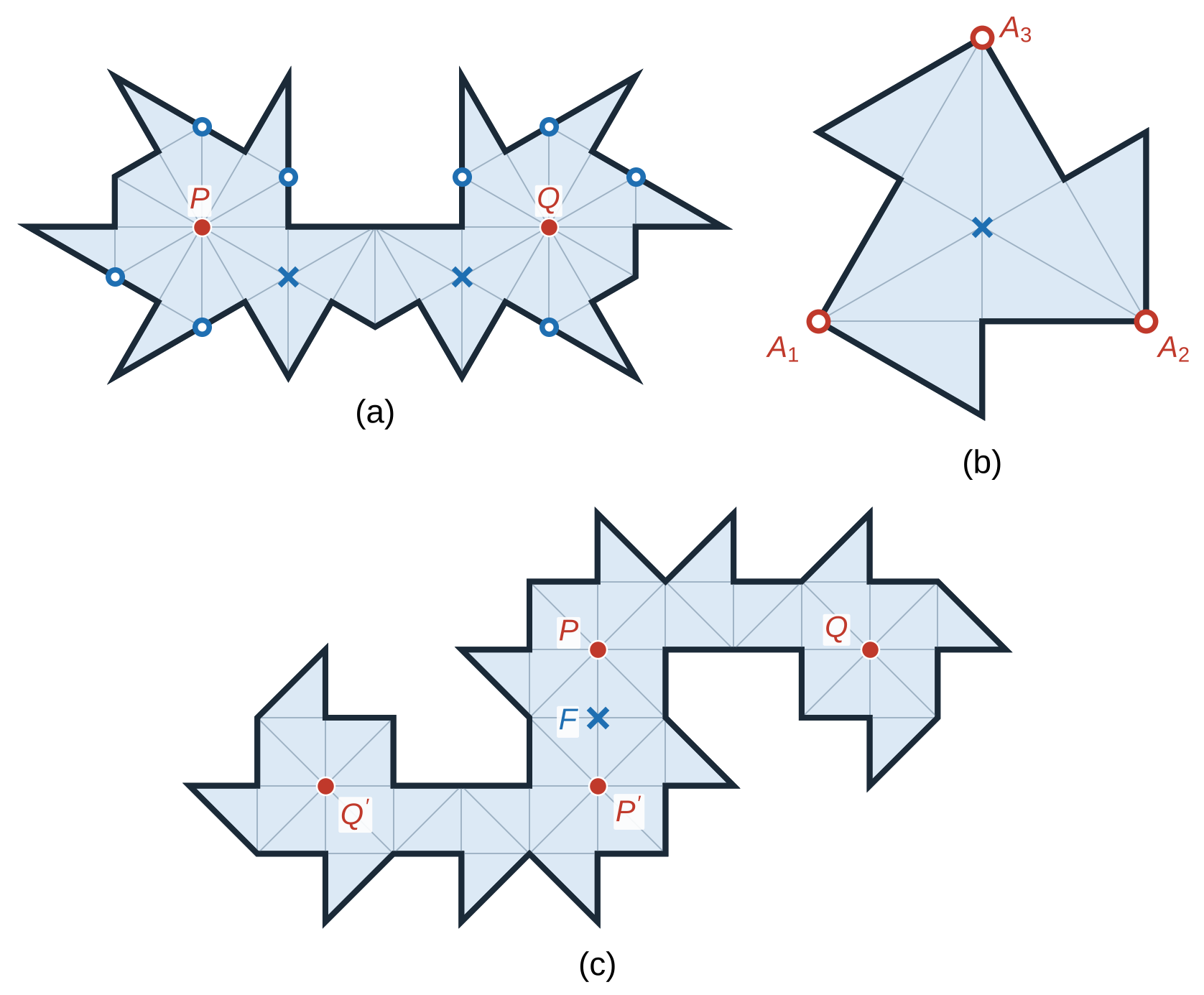}
\caption{Rooms in which some tile vertices are not corners. (a)~The 27-gon of Proposition~\ref{prop:27}, with mutually invisible interior points $P$ and $Q$; crosses mark its two interior $B$-points and small circles its eight flat $B$-points. (b)~The nonagon of Proposition~\ref{prop:nonagon}, with its three invisible corners $A_1$, $A_2$, $A_3$; the cross marks the interior $B$-point. (c)~The 42-gon of Example~\ref{ex:42}; no light ray joins a point of $\{P,P'\}$ to a point of $\{Q,Q'\}$, and the cross marks the centre $F$ of the rotation.}
\label{fig:relax}
\end{figure}

\section{A ring of trapezoids}\label{sec:trapezoid}

The proof of Proposition~\ref{prop:star}(a) used folding. Folding works for a room made of copies of one polygon, the \emph{base tile}, in which any two copies with a common side are mirror images in it: each copy is mapped onto the base tile by the reflections in the sides crossed along a chain of copies, and its vertices receive the labels of the corresponding vertices of the base tile. This needs the labels to be consistent: each vertex must receive the same label along every chain. A hole can break this consistency, and then a light ray can return to a corner whose image in the base tile is invisible.

\begin{figure}[htbp]
\centering
\includegraphics[width=0.92\textwidth]{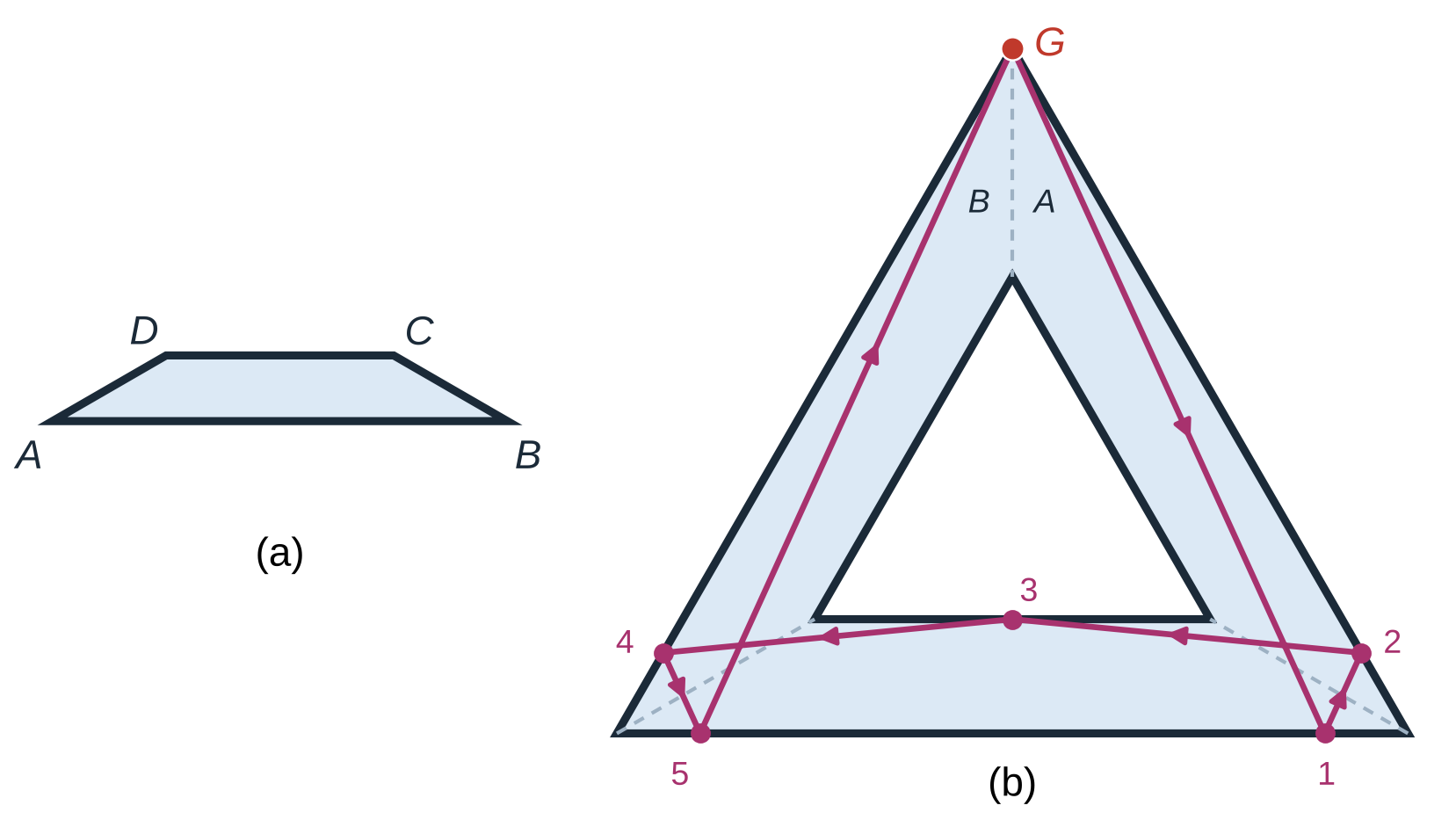}
\caption{(a)~The trapezoid $\tau$, with angles $30\dg$ at $A$ and $B$. (b)~Three copies of $\tau$ around a triangular hole, and the light ray of Proposition~\ref{prop:return}, which leaves the corner $G$ and returns to it after five reflections. Dashed lines are the edges between the copies. The two copies at $G$ carry the labels $A$ and $B$ there.}
\label{fig:trapezoid}
\end{figure}

Let $\tau$ be the isosceles trapezoid $ABCD$ with long base $AB$, short base $CD$, angles $30\dg$ at $A$ and $B$, and angles $150\dg$ at $C$ and $D$ (Figure~\ref{fig:trapezoid}(a)). The angle at $A$ is $30\dg$, which divides $90\dg$, and every angle of $\tau$ is a multiple of $30\dg$; by Theorem~1 of Tokarsky and Boese~\cite{TB}, no light ray in $\tau$ starts and ends at $A$. Reflect $\tau$ in its two legs. The three trapezoids fit together around an equilateral triangular hole bounded by their short bases, and their long bases form an equilateral triangle: at each outer corner two angles of $30\dg$ meet, and at each corner of the hole two angles of $150\dg$ leave $60\dg$. Label the vertices of the three trapezoids by reflection from $\tau$. Going once around the hole composes three reflections, which map $\tau$ onto itself reversing orientation, so the composition is the symmetry of $\tau$ that exchanges $A$ with $B$ and $C$ with $D$. A light ray that leaves a corner of the ring and returns to it after going once around the hole therefore folds into a light ray in $\tau$ from $A$ to $B$, which Theorem~1 of~\cite{TB} does not exclude. The ring and this observation are in the second author's letter~\cite[Figure~3(d),(f)]{Tok96}, which also draws a light ray in the ring that returns to a corner.

\begin{proposition}\label{prop:return}
In coordinates in which $(x,y)$ denotes the point $(x,\sqrt3\,y)$, let $\tau$ have vertices $A=(-6,0)$, $B=(6,0)$, $C=(3,1)$, $D=(-3,1)$. The ring is the triangle with vertices $(-6,0)$, $(6,0)$, $(0,6)$ with the open triangle with vertices $(-3,1)$, $(3,1)$, $(0,4)$ removed. The polygonal path from $G=(0,6)$ through
\[
\bigl(\tfrac{90}{19},0\bigr),\quad\bigl(\tfrac{90}{17},\tfrac{12}{17}\bigr),\quad(0,1),\quad\bigl(-\tfrac{90}{17},\tfrac{12}{17}\bigr),\quad\bigl(-\tfrac{90}{19},0\bigr)
\]
back to $G$ is a light ray with five reflections that goes once around the hole. The ray leaves $G$ inside the copy of $\tau$ in which $G$ is the vertex $A$, and it returns to $G$ inside the copy in which $G$ is the vertex $B$.
\end{proposition}

\begin{proof}
A direct computation with rational arithmetic, in the metric $dx^2+3\,dy^2$, verifies the law of reflection at each of the five points, which lie in the relative interiors of sides, and verifies that each segment lies in the ring and meets no corner.
\end{proof}

The outer boundary of the ring can be changed: rooms made of four and of five copies of $\tau$, with a pentagonal and a heptagonal outer boundary around the same triangular hole, contain similar returns to a corner with $5$ and $19$ reflections.

In rooms made of copies of $\tau$, a single hole exchanges labels only in this way.

\begin{proposition}\label{prop:onehole}
Let $R$ be a union of finitely many copies of $\tau$ with pairwise disjoint interiors and connected interior, such that any two copies that meet in a segment of positive length share a whole side and are mirror images of each other in it. Suppose that $R$ has exactly one hole, whose boundary is disjoint from the outer boundary. If carrying the labels of a copy once around the hole exchanges $A$ and $B$, then the hole is an equilateral triangle bounded by three short bases, as in Figure~\textup{\ref{fig:trapezoid}(b)}.
\end{proposition}

\begin{proof}
A copy has at most four neighbours, along its long base, its short base and its two legs. A neighbour along the short base excludes a neighbour along either leg, since three angles of $150\dg$ would meet at an endpoint of the short base. If a copy has neighbours along both legs, these neighbours share a leg, and the three copies form the ring of Figure~\ref{fig:trapezoid}(b); its triangular hole cannot be filled, since a copy along a short base is excluded, so it is the hole of $R$. Suppose now that no copy has neighbours along both legs. Then every copy has at most two neighbours, and a copy with two neighbours has one along its long base and one along its short base or a leg. In any cycle of the adjacency graph of the copies, edges along long bases therefore alternate with edges of the other kinds, so every cycle has even length and every closed walk has even length. Carrying the labels around the hole composes an even number of reflections, which maps $\tau$ onto itself and keeps its orientation. The only such symmetry of $\tau$ is the identity, so the labels return unchanged, contrary to the hypothesis.
\end{proof}

\section{Remarks on the computations}\label{sec:computations}

The searches and verifications use exact integer or rational arithmetic. Two programs that construct examples use floating-point methods, an integer linear program for the assemblies of Section~\ref{sec:many} and a numerical ray search for the trapezoid ring, and exact checks verify their output. The searches behind Section~\ref{sec:pinwheel} and Theorems~\ref{thm:posts}, \ref{thm:hole} and~\ref{thm:three} are finite by the bounds given in their proofs and use the search method of Theorem~\ref{thm:tri}. The searches behind Theorems~\ref{thm:post}, \ref{thm:pinched}, \ref{thm:tri}, \ref{thm:sq} and~\ref{thm:min27} are complete because of the bounds proved above. Theorems~\ref{thm:post} and~\ref{thm:pinched} reduce to $64+576$ and $256$ explicit candidates, Theorem~\ref{thm:tri} to $77$ positions of the second source with at most $28$ tiles, and Theorem~\ref{thm:sq} to $35$ positions with at most $21$ squares, searched with up to $22$ by one program and up to $21$ by a second. Theorem~\ref{thm:min27} reduces to $13$ configurations, which an exhaustive branching program rejects after visiting $1{,}527$ branch nodes and an independent program using the OR-Tools constraint solver also rejects. For $27$ sides the branching program searches $103$ configurations and lists $50$ rooms. For every room in our lists a separate program checks, with integer orientation tests, that the boundary is simple, that the corner condition is satisfied, and that the number of sides is as stated; the types are the orbits of the stated groups. For the examples of Sections~\ref{sec:many}--\ref{sec:trapezoid} it checks the coordinates, the module structure, the reflections in the sides and the light rays. The checks support the proofs of the finite cases; the general lemmas are proved in the text.

The programs and data accompany this paper. The command \texttt{python run\_checks.py} runs the default checks in under a minute, using the Python standard library and Node.js; the OR-Tools program is optional. The command \texttt{python run\_checks.py --full-census} also repeats the longer searches: those of Sections~\ref{sec:pinwheel} and~\ref{sec:holefree} and of Theorems~\ref{thm:posts}, \ref{thm:hole} and~\ref{thm:three}. Table~\ref{tab:code} lists the main programs.

\begin{table}[htbp]
\centering\small
\begin{tabular}{@{}>{\raggedright\arraybackslash}p{0.30\textwidth}>{\raggedright\arraybackslash}p{0.64\textwidth}@{}}
\toprule
Result & Programs\\
\midrule
Section~\ref{sec:one} & \nolinkurl{code/single_source/verify_single_source.py}, \nolinkurl{one_source_census.py}, \nolinkurl{code/consistency_checks.py}\\
Theorem~\ref{thm:post}, Corollary~\ref{cor:partition} & \nolinkurl{code/rooms_with_holes/enumerate_posts.py}, \nolinkurl{verify_posts.py}\\
Theorem~\ref{thm:posts} & \nolinkurl{code/rooms_with_holes/post_search.py}\\
Theorem~\ref{thm:pinched} & \nolinkurl{code/rooms_with_holes/enumerate_pinched.py}\\
Theorem~\ref{thm:hole} & \nolinkurl{code/rooms_with_holes/hole_search.py}\\
Remark~\ref{rem:wolecki} & \nolinkurl{code/rooms_with_holes/wolecki22_match.py}\\
Theorem~\ref{thm:tri} & \nolinkurl{code/hole_free_census/room_search_shared_grid.py}, \nolinkurl{validate_census.py}\\
Theorem~\ref{thm:sq} & \nolinkurl{code/hole_free_census/room_search.py}, \nolinkurl{enumerate_orthogonal.py}, \nolinkurl{validate_census.py}\\
Section~\ref{sec:many} & \nolinkurl{code/assemblies/verify_snakes.py}, \nolinkurl{verify_assemblies.py}\\
Theorem~\ref{thm:three} & \nolinkurl{code/hole_free_census/three_sources.py}\\
Proposition~\ref{prop:27}, Theorem~\ref{thm:min27} & \nolinkurl{code/minimum27/example27/verify_example.py}, \nolinkurl{code/minimum27/verify_minimum.js}, \nolinkurl{independent_cp_check.py}, \nolinkurl{enumerate27.js}, \nolinkurl{check_27gons.py}\\
Proposition~\ref{prop:nonagon} & \nolinkurl{code/nonagon/geometry.py}\\
Example~\ref{ex:42} & \nolinkurl{code/symmetry_rooms/verify_symmetry_room.py}\\
Section~\ref{sec:trapezoid} & \nolinkurl{code/trapezoid_rings/verify_examples.py}\\
Coordinates of Appendix~\ref{app:plates} & \nolinkurl{data/hole_free_catalogue.json}\\
Figures & \nolinkurl{code/figures/make_paper_figures.py}\\
\bottomrule
\end{tabular}
\caption{Programs that carry out or check the computations.}
\label{tab:code}
\end{table}

\section{Concluding remarks}\label{sec:conclusion}

The corner condition turns questions about light into questions about counting tiles. Two counting lemmas, Lemmas~\ref{lem:euler} and~\ref{lem:wings}, bound the number of tiles of a lattice room by its number of sides, and this bound forces the sources to lie close together. This is why the pinwheel, the cross $\Sigma$ and the two-source rooms of Sections~\ref{sec:devices} and~\ref{sec:holefree} can be listed completely, and why, in the rooms listed there, the two sources are near each other: at distance $2$ with posts or a partition, $2\sqrt2$ with a hole, $4$ or $2\sqrt5$ in the hole-free lattice rooms on $\T$, and $4$, $2\sqrt5$ or $4\sqrt2$ on $\Sq$. The post, the partition and the pinched hole each save sides in their own way. A post absorbs light at the $C$-point shared by two sources at distance $2$; a partition is such a post joined to the outer boundary; and a pinched hole lets the tiles of the two sources meet along a diagonal. Up to three posts, each further post saves one side, and a hole that does not touch the outer boundary costs four sides more than a pinched one.

The same arithmetic explains the rooms of Section~\ref{sec:relax}. Lemma~\ref{lem:blocking} asks for corners only at the points of $\Lambda\setminus2\Lambda$, and the tiling decides which points these are. On the tiling by triangles with angles $30\dg$, $60\dg$ and $90\dg$ they are the $C$-points alone: the 27-gon keeps two $B$-points in its interior and eight on its sides, and in the nonagon three $C$-points suffice to make its three corners of angle $90\dg$ mutually invisible. In the 42-gon a rotation through $180\dg$ takes the place of the one missing corner. The trapezoid ring marks a limit of the folding method: around a hole the labels of the tiles can change, and a light ray can then return to a corner that is invisible in the base tile.

We close with some questions.
\begin{enumerate}[label=\arabic*.]
\item The pinwheel is the only hole-free lattice room on $\T$ with one source and $12$ sides, and none of the three such rooms with $13$ sides contains it. Section~\ref{sec:pinwheel} gives the number of types up to $20$ sides. How does this number grow?
\item With at most four posts, a room on $\T$ with two sources has at least $16$ sides (Theorem~\ref{thm:posts}). Do more posts lower this number, and what is the least number of sides with any number of posts?
\item Theorems~\ref{thm:pinwheel}, \ref{thm:tri} and~\ref{thm:three} show that hole-free lattice rooms on $\T$ with $n=1$, $2$ and $3$ sources have at least $12n$ sides, and the pinwheel and the snakes of Proposition~\ref{prop:snakes} have exactly $12n$ sides. Is $12n$ the least number for every $n$? Which trees, and which other planar graphs, are graphs of assemblies?
\item Tokarsky~\cite[Example~3]{Tok95} asked how few sides a polygon with two mutually invisible interior points can have. Wolecki's room~\cite{Wolecki} and Theorem~\ref{thm:pinched} give $22$ for rooms with a pinched hole, and Thomson's room~\cite{Thomson} has $21$ if sides are counted as maximal straight segments (Remark~\ref{rem:pinchcount}). For a simple polygon, Theorem~\ref{thm:tri} gives $24$ within lattice rooms on $\T$, and Theorem~\ref{thm:min27} gives $27$ within its class on $\K$. Can a simple polygon with fewer than $24$ sides have two mutually invisible interior points?
\item Proposition~\ref{prop:onehole} shows that one hole exchanges labels only in the ring of three trapezoids. Which arrangements of several holes exchange labels? Can a light ray in such a room join two interior points whose images in the base tile are mutually invisible?
\end{enumerate}

\section*{Acknowledgements}
The proofs were generated and the draft of this paper was written with GPT-6 (OpenAI), and the draft was checked and rewritten with Claude Opus 5.5 (Anthropic). This work was done under the guidance of the authors from March 2026 onwards, and earlier versions of these models were used as well.

\appendix

\clearpage
\section{Catalogues}\label{app:plates}

\begin{figure}[!htbp]
\centering
\includegraphics[width=0.9\textwidth]{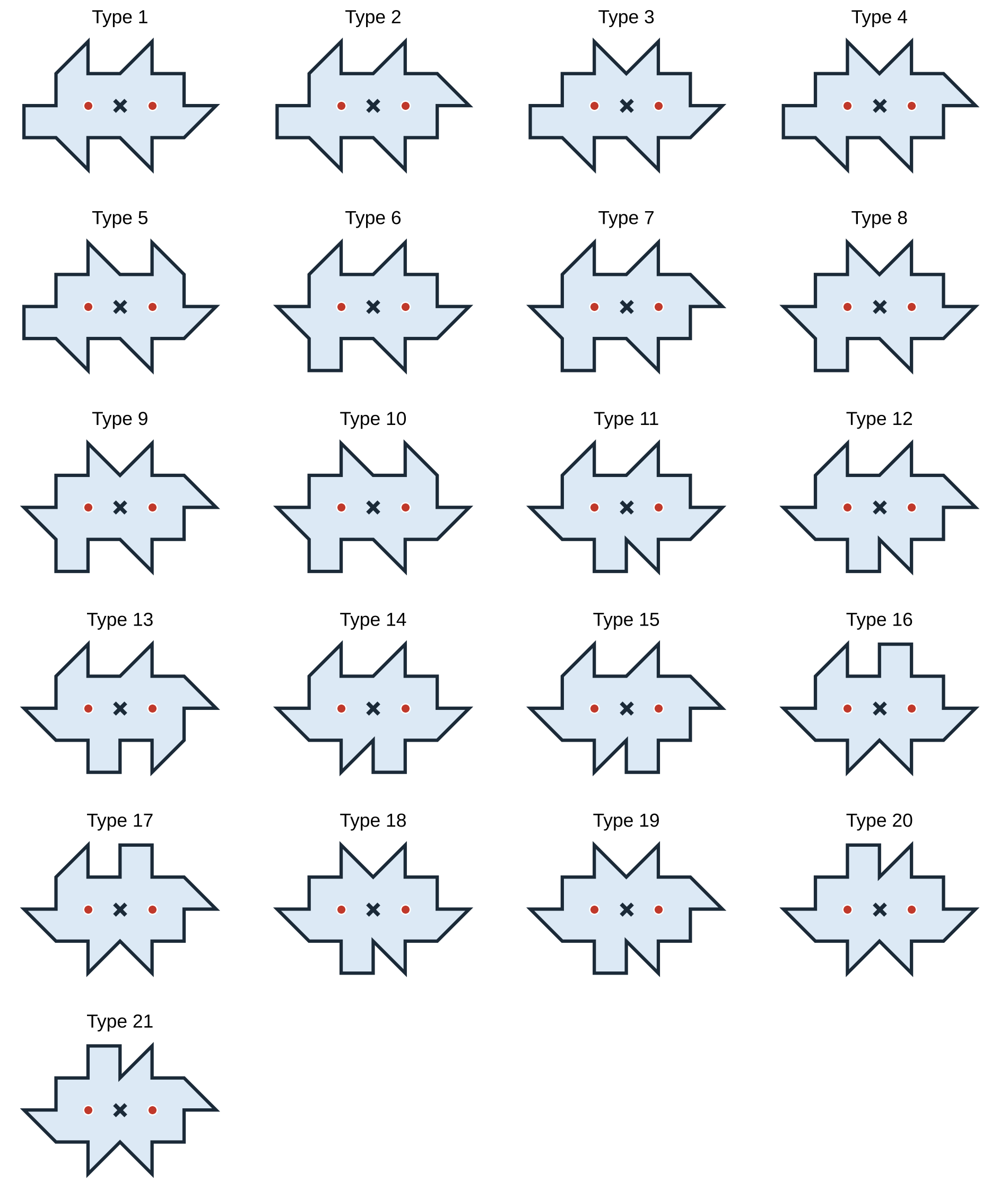}
\caption{The twenty-one post rooms with 19 sides (Theorem~\ref{thm:post}).}
\label{fig:post19}
\end{figure}

\begin{figure}[p]
\centering
\includegraphics[width=\textwidth]{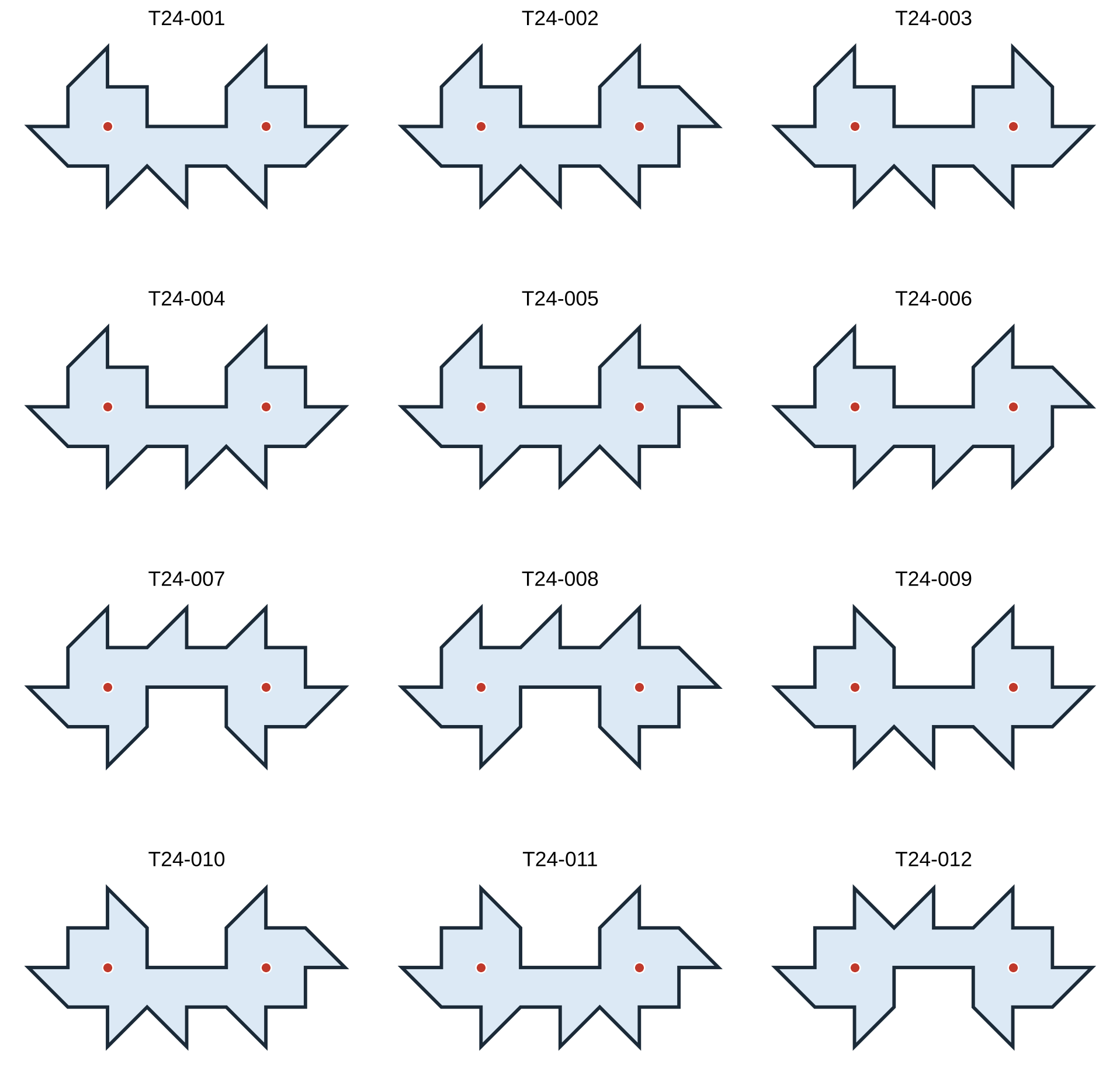}
\caption{The twelve hole-free lattice rooms on $\T$ with 24 sides, T24-001 to T24-012 (Theorem~\ref{thm:tri}). The sources are $(0,0)$ and $(4,0)$. T24-003, T24-008 and T24-012 are the rooms of Castro~\cite{Castro}, of the second author~\cite[Figure~2]{Tok96} and of Thomson~\cite{Thomson}.}
\label{fig:t24}
\end{figure}

\begin{figure}[p]
\centering
\includegraphics[width=\textwidth]{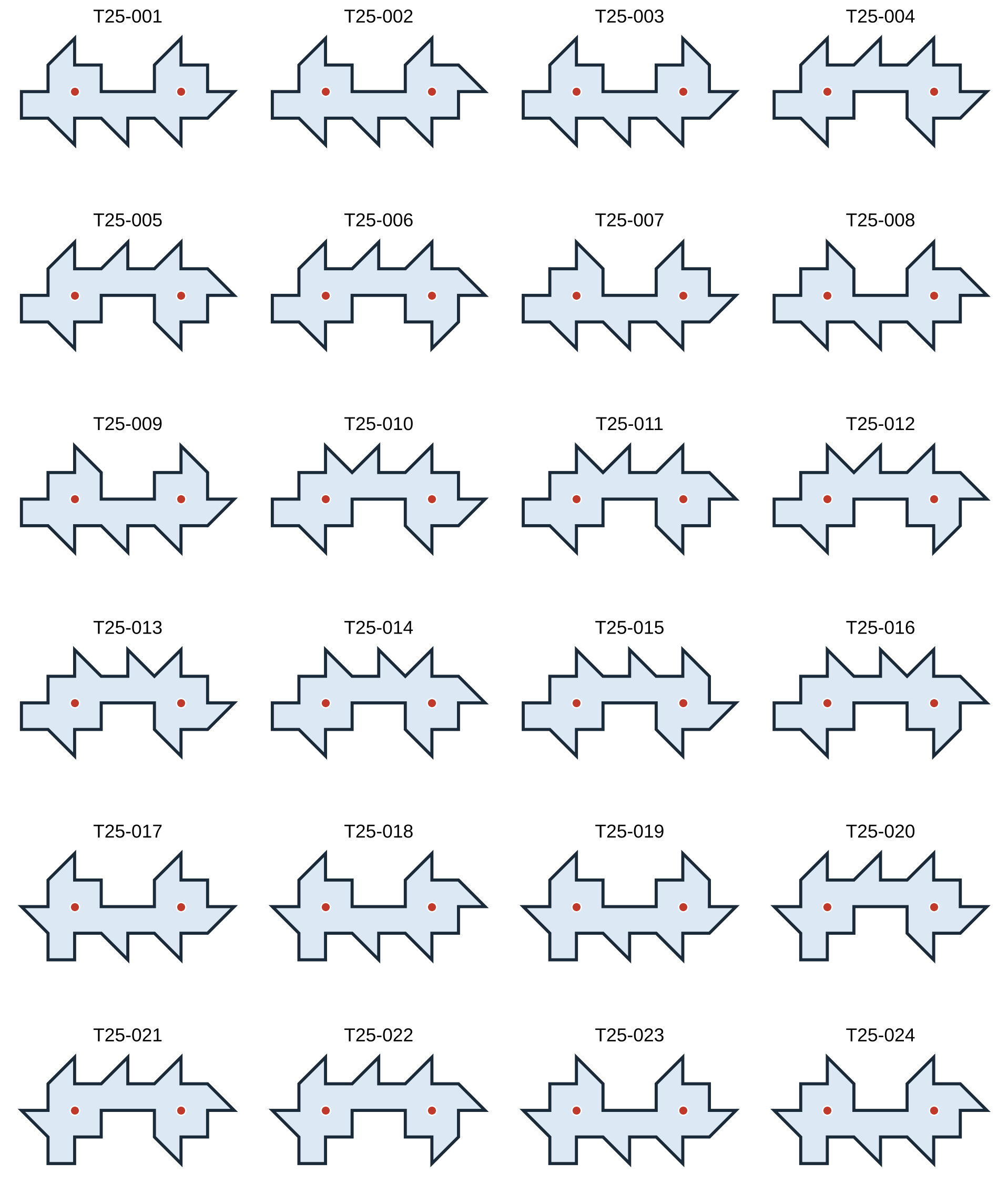}
\caption{Hole-free lattice rooms on $\T$ with 25 sides and sources $(0,0)$, $(4,0)$: T25-001 to T25-024.}
\label{fig:t25a}
\end{figure}

\begin{figure}[p]
\centering
\includegraphics[width=\textwidth]{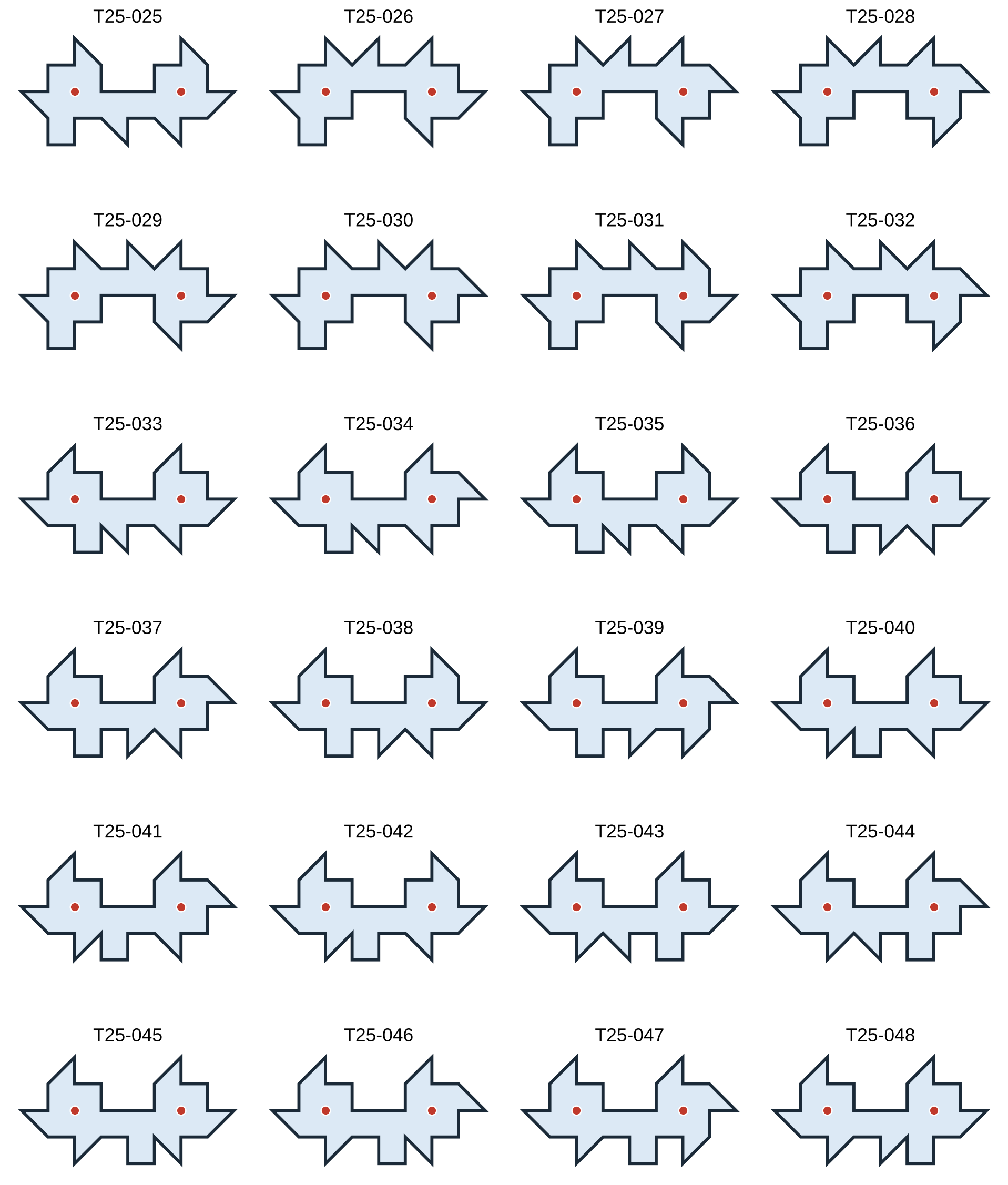}
\caption{Hole-free lattice rooms on $\T$ with 25 sides and sources $(0,0)$, $(4,0)$: T25-025 to T25-048.}
\label{fig:t25b}
\end{figure}

\begin{figure}[p]
\centering
\includegraphics[width=\textwidth]{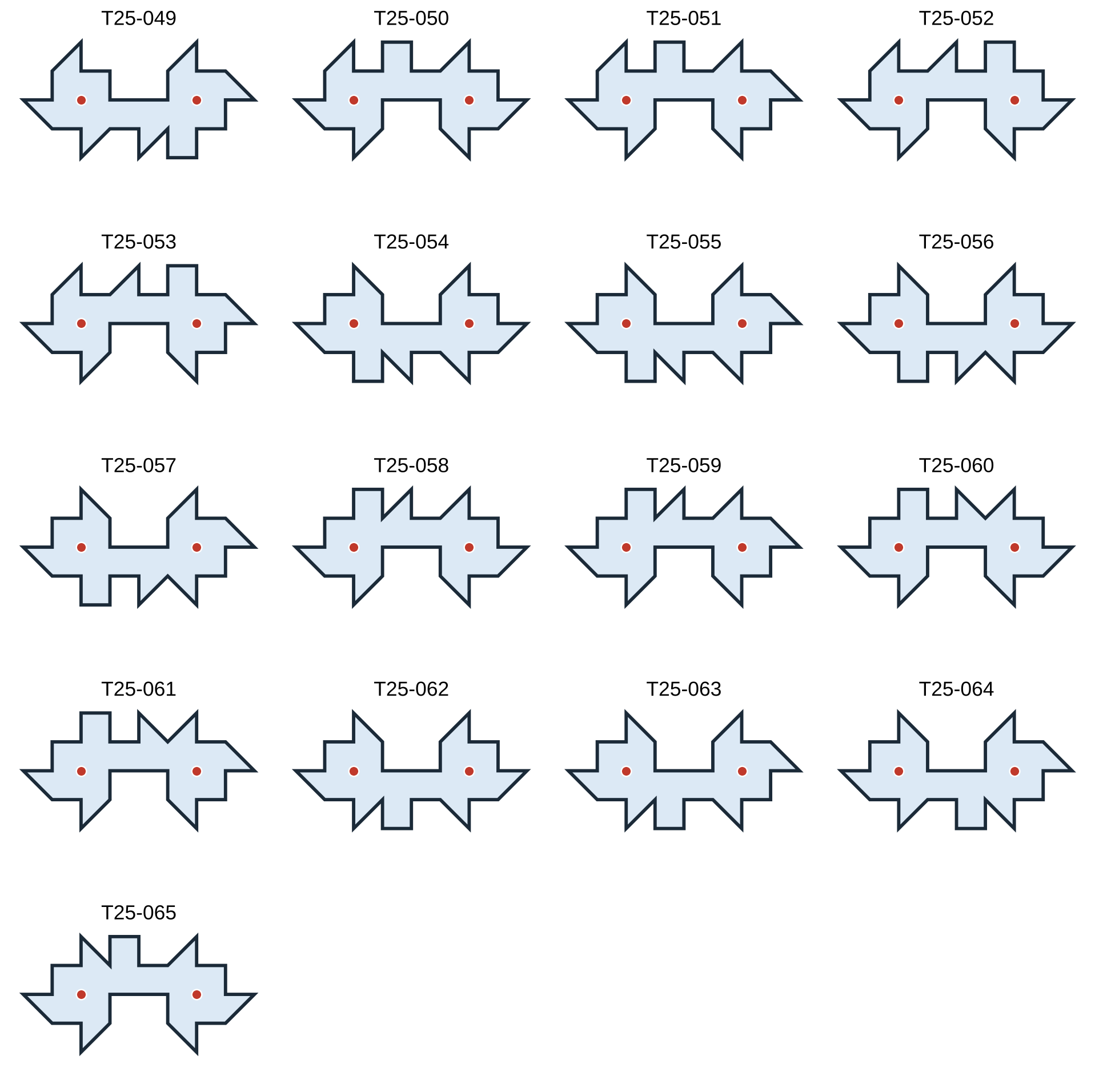}
\caption{Hole-free lattice rooms on $\T$ with 25 sides and sources $(0,0)$, $(4,0)$: T25-049 to T25-065.}
\label{fig:t25c}
\end{figure}

\begin{figure}[p]
\centering
\includegraphics[width=\textwidth]{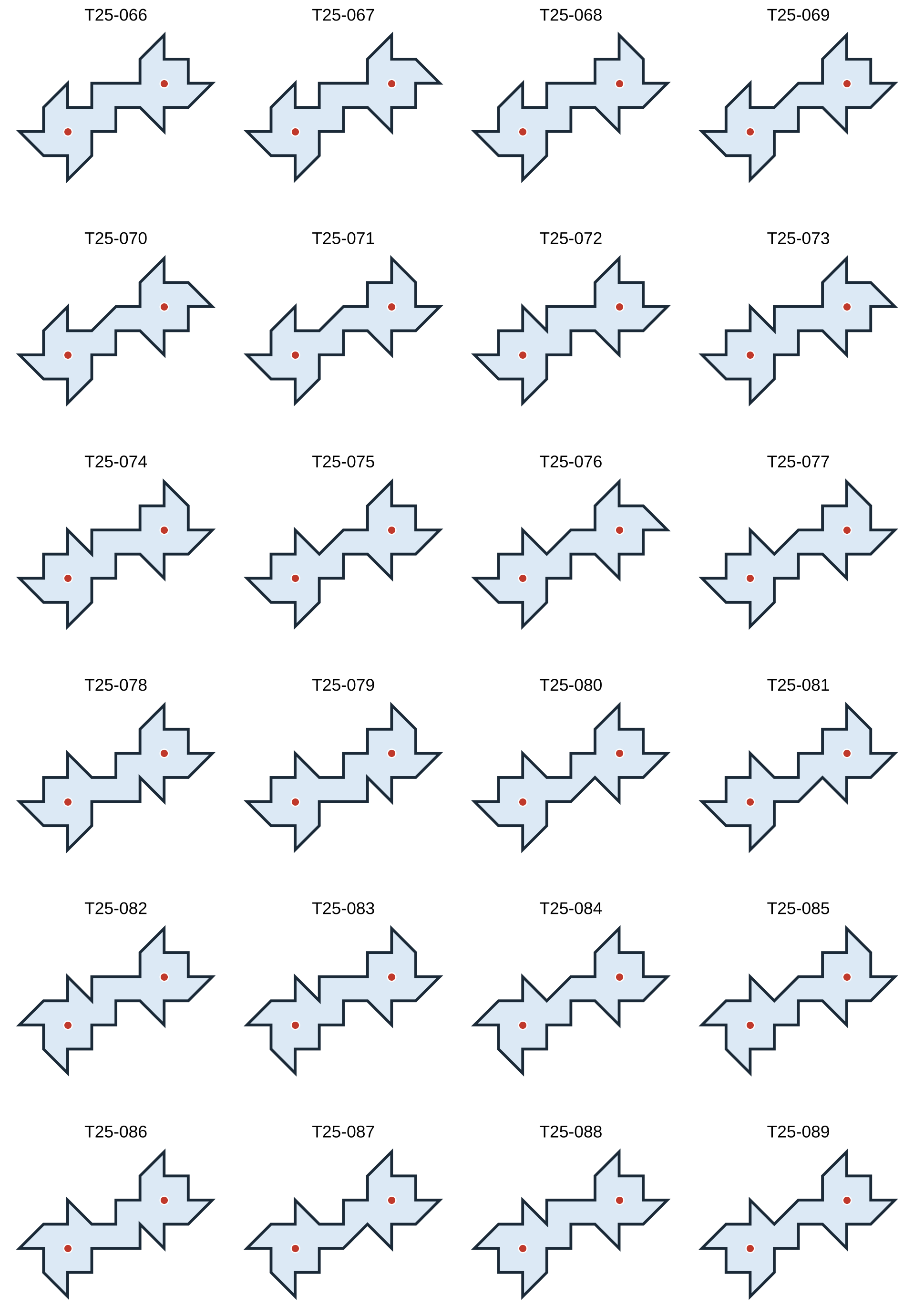}
\caption{Hole-free lattice rooms on $\T$ with 25 sides and sources $(0,0)$, $(4,2)$: T25-066 to T25-089.}
\label{fig:t25o}
\end{figure}

\begin{figure}[p]
\centering
\includegraphics[width=\textwidth]{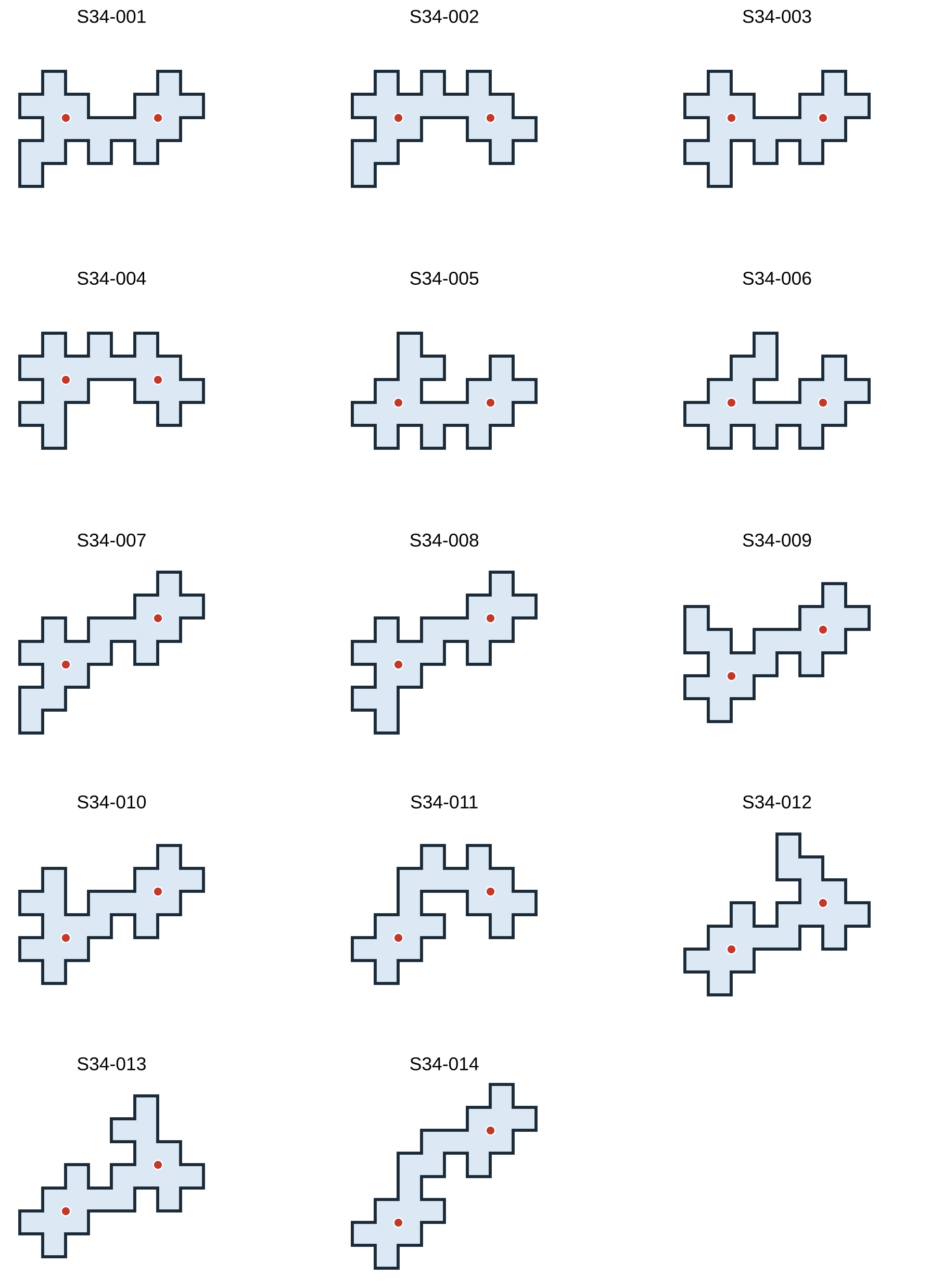}
\caption{The fourteen hole-free lattice rooms on $\Sq$ with 34 sides (Theorem~\ref{thm:sq}). The sources are $(0,0)$ and $(4,0)$ for S34-001 to S34-006, $(0,0)$ and $(4,2)$ for S34-007 to S34-013, and $(0,0)$ and $(4,4)$ for S34-014.}
\label{fig:s34}
\end{figure}

\begin{figure}[p]
\centering
\includegraphics[width=0.78\textwidth]{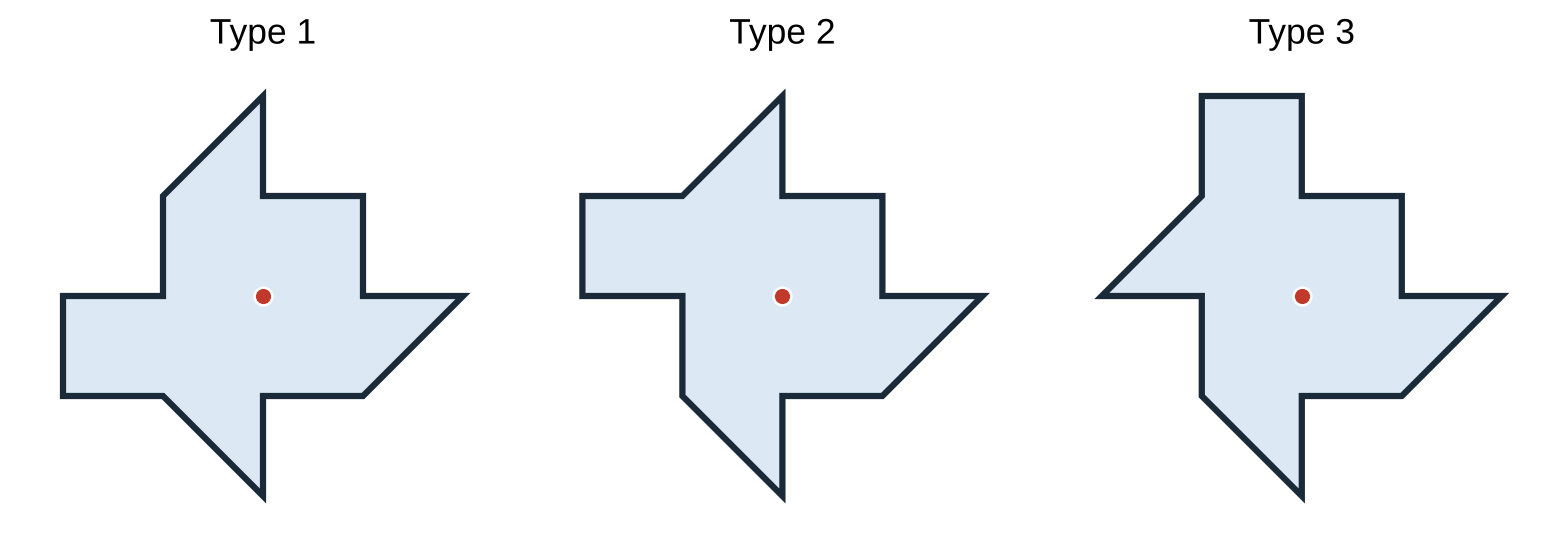}
\caption{The three hole-free lattice rooms on $\T$ with one source and 13 sides (Section~\ref{sec:pinwheel}). None of them contains the pinwheel.}
\label{fig:onesource13}
\end{figure}

\begin{figure}[p]
\centering
\includegraphics[width=0.78\textwidth]{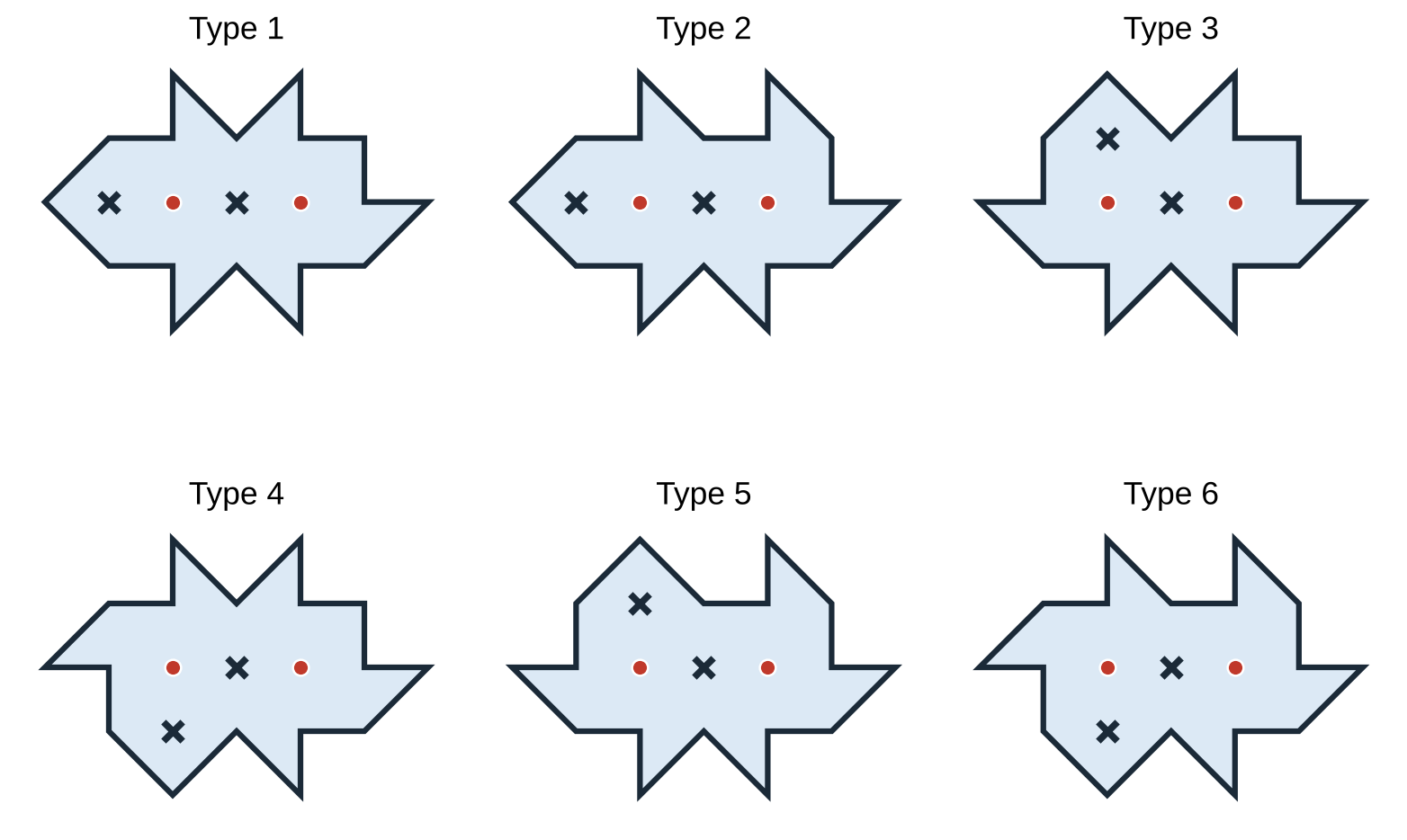}
\caption{The six rooms with two sources, two posts and 17 sides (Theorem~\ref{thm:posts}). Crosses mark the posts; the sources are $(0,0)$ and $(2,0)$.}
\label{fig:posts2}
\end{figure}

\begin{figure}[p]
\centering
\includegraphics[width=\textwidth]{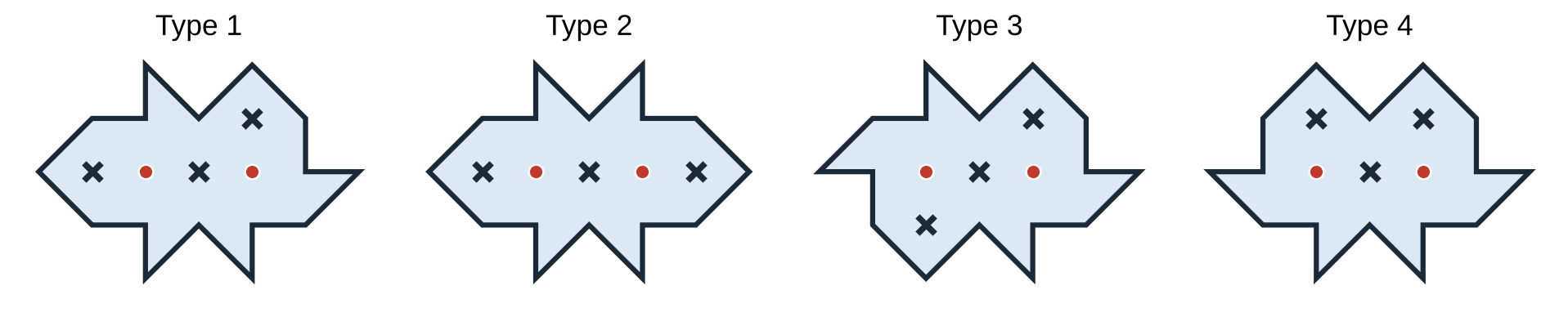}
\caption{The four rooms with two sources, three posts and 16 sides (Theorem~\ref{thm:posts}). Crosses mark the posts.}
\label{fig:posts3}
\end{figure}

\begin{figure}[p]
\centering
\includegraphics[width=\textwidth]{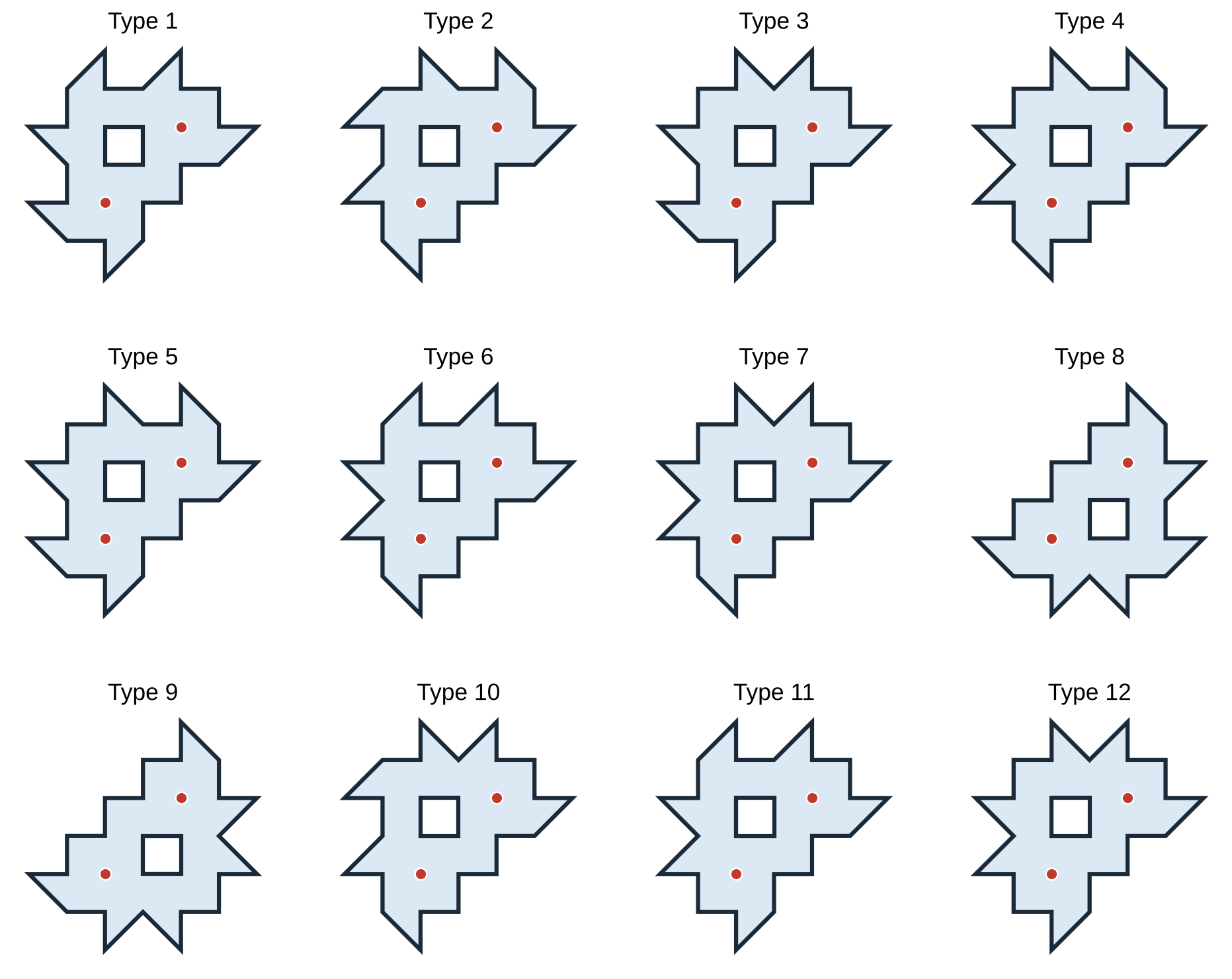}
\caption{The twelve lattice rooms on $\T$ with two sources, one hole and 26 sides (Theorem~\ref{thm:hole}).}
\label{fig:hole}
\end{figure}

\begin{figure}[p]
\centering
\includegraphics[width=\textwidth]{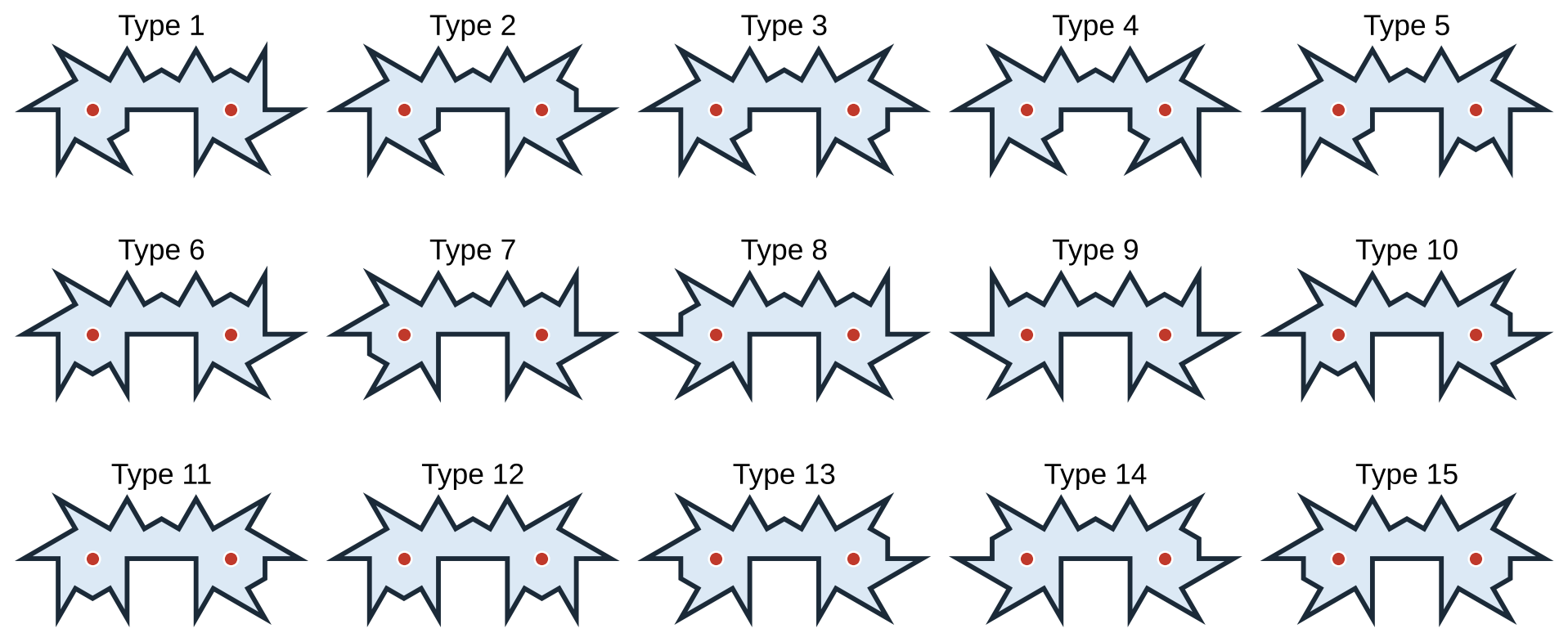}
\caption{The fifteen rooms on $\K$ with two interior $A$-points and 27 sides in which every $C$-point is a corner (Theorem~\ref{thm:min27}). The sources are $P$ and $Q$ of Proposition~\ref{prop:27}.}
\label{fig:k27}
\end{figure}

Figures~\ref{fig:post19}--\ref{fig:k27} show one room of each type of Theorems~\ref{thm:post}, \ref{thm:tri} and~\ref{thm:sq}, of Section~\ref{sec:pinwheel} for $13$ sides, and of Theorems~\ref{thm:posts}, \ref{thm:hole} and~\ref{thm:min27}; the smaller catalogues are Figures~\ref{fig:post18}, \ref{fig:partition} and~\ref{fig:pinched}. Dots mark the sources. Within a figure all rooms are drawn to the same scale. The exact coordinates and tiles of all rooms, with the number of rooms of each type, are in the files that accompany the paper. The labels T24, T25 and S34 follow a fixed order: by the position of the second source, and then by the sorted list of tiles of a canonical representative.

\clearpage
\section{The straight and bent snakes}\label{app:snakes}

We write $\{a,\dots,b\}$ for the set of integers from $a$ to $b$, and $a,a+2,\dots,b$ or $a,a+4,\dots,b$ for arithmetic progressions that include both ends.

\subsection{The straight snake \texorpdfstring{$\Gamma_n$}{Gamma n}}
Let $n\ge2$ and $d=4(n-1)$. The boundary of $\Gamma_n$ consists of an upper path from $(-2,0)$ to $(d+2,0)$ and a lower path back. The upper path passes through
\[
(-2,0),\ (-1,1),\ (0,1),\ (0,2),\ (1,1),\ (2,1),\ (2,2),\ (3,1),\ (4,2),
\]
then, for $k=2,3,\dots,2n-3$, through $(2k,1)$, $(2k+1,1)$, $(2k+2,2)$, and ends with $(d,1)$, $(d+1,1)$, $(d+2,0)$. The lower path passes through
\[
(d+1,0),\ (d+1,-1),\ (d,-2),\ (d,-1),
\]
then, for $c=d-4,d-8,\dots,4$, through $(c+3,-1)$, $(c+3,0)$, $(c+1,0)$, $(c+1,-1)$, $(c,-2)$, $(c,-1)$, and ends with
\[
(3,-1),\ (3,0),\ (1,0),\ (1,-1),\ (0,-1),\ (0,-2),\ (-1,-1),\ (-1,0).
\]
For $n=2$ the two repeated blocks are empty. The upper path has $9+3(2n-4)+3=6n$ vertices and the lower path $4+6(n-2)+8=6n$, so the boundary has $12n$ vertices.

Every edge lies on a mirror line~\eqref{eq:mirrors}, and consecutive edges are not collinear, so all $12n$ vertices are corners. Along the upper path, $x$ increases strictly on every edge that is not vertical, and each vertical line contains at most one vertical edge; the same is true for the lower path with $x$ decreasing. Hence each path is simple. The upper path has $y>0$ except at its endpoints, and the lower path has $y\le0$, so the two paths meet only where the boundary closes, and $\Gamma_n$ is a simple polygon. Reading off the two paths, the lattice points of $\Gamma_n$ are those listed in Table~\ref{tab:straightrows}.

\begin{table}[htbp]
\centering\small
\begin{tabular}{@{}cl@{}}
\toprule
$y$ & $x$-coordinates of the lattice points of $\Gamma_n$\\
\midrule
$2$ & $0,2,\dots,d$\\
$1$ & $\{-1,\dots,d+1\}$\\
$0$ & $\{-2,\dots,d+2\}$\\
$-1$ & $-1,1,\dots,d+1$ and $0,4,\dots,d$\\
$-2$ & $0,4,\dots,d$\\
\bottomrule
\end{tabular}
\caption{Lattice points of the straight snake $\Gamma_n$, with $d=4(n-1)$.}
\label{tab:straightrows}
\end{table}

The interior lattice points are exactly $(4j,0)$ for $0\le j<n$. All other points in the table lie on the boundary, and the boundary lattice points that are not vertices are the $A$-points $(4j+2,0)$, $0\le j<n-1$. So $\Gamma_n$ is a lattice room with $M=n$, $f=n-1$ and $N=12n$, and~\eqref{eq:euler} gives $T=15n-3$. The parities of the vertices give $3n+1$ $A$-points, $4n$ $B$-points and $5n-1$ $C$-points.

\subsection{The bent snake \texorpdfstring{$\Gamma'_n$}{Gamma' n}}
Let $n\ge3$ and $d=4(n-2)$. The boundary of $\Gamma'_n$ consists of three paths. The bend runs from $(-2,0)$ to $(1,0)$ through the $18$ vertices
\[
\begin{gathered}
(-2,0),\ (-1,1),\ (-2,2),\ (-1,2),\ (-1,3),\ (-2,4),\ (-1,4),\ (-1,5),\ (0,6),\\
(0,5),\ (1,5),\ (1,4),\ (2,4),\ (1,3),\ (0,3),\ (0,1),\ (1,1),\ (1,0).
\end{gathered}
\]
The upper path passes, for $k=1,\dots,n-2$ and $c=4k$, through $(c-1,0)$, $(c-1,1)$, $(c,2)$, $(c,1)$, $(c+1,1)$, $(c+1,0)$, and ends at $(d+2,0)$; it has $6(n-2)+1$ vertices. The lower path passes, for $e=d,d-2,\dots,4$, through $(e+1,-1)$, $(e,-1)$, $(e,-2)$, and ends with
\[
(3,-1),\ (2,-2),\ (2,-1),\ (1,-1),\ (0,-2),\ (0,-1),\ (-1,-1),\ (-1,0);
\]
it has $3(2n-5)+8$ vertices. In all there are $18+6(n-2)+1+3(2n-5)+8=12n$ vertices.

As before, every edge lies on a mirror line and consecutive edges are not collinear. The bend lies in the strip $-2\le x\le2$ with $y\ge0$, meets the line $y=0$ only at its endpoints, and is simple, as one checks directly. The upper path lies in the half-plane $x\ge3$ with $y\ge0$, and the lower path has $y<0$ except at its last vertex $(-1,0)$. Along the upper path $x$ never decreases, along the lower path $x$ never increases, and on each path every vertical line contains at most one vertical edge, so each path is simple, and no two of the three paths meet. The three remaining edges, from $(1,0)$ to $(3,0)$, from $(d+2,0)$ to $(d+1,-1)$ and from $(-1,0)$ to $(-2,0)$, join the paths, and each meets the rest of the boundary only at its endpoints. Hence $\Gamma'_n$ is a simple polygon. Its lattice points are listed in Table~\ref{tab:bentrows}.

\begin{table}[htbp]
\centering\small
\begin{tabular}{@{}cl@{}}
\toprule
$y$ & $x$-coordinates of the lattice points of $\Gamma'_n$\\
\midrule
$6$ & $0$\\
$5$ & $-1,0,1$\\
$4$ & $\{-2,\dots,2\}$\\
$3$ & $-1,0,1$\\
$2$ & $-2,-1,0$ and $4k$ for $1\le k\le n-2$\\
$1$ & $-1,0,1$ and $4k-1,4k,4k+1$ for $1\le k\le n-2$\\
$0$ & $\{-2,\dots,d+2\}$\\
$-1$ & $\{-1,\dots,d+1\}$\\
$-2$ & $0,2,\dots,d$\\
\bottomrule
\end{tabular}
\caption{Lattice points of the bent snake $\Gamma'_n$, with $d=4(n-2)$.}
\label{tab:bentrows}
\end{table}

The interior lattice points are exactly $(4j,0)$ for $0\le j\le n-2$ and $(0,4)$. All other points in the table lie on the boundary, and the boundary lattice points that are not vertices are the $A$-points $(0,2)$ and $(4k-2,0)$, $1\le k\le n-2$; there are $n-1$ of them. This proves the corner condition, $M=n$ and $f=n-1$, and the parities of the vertices and~\eqref{eq:euler} give the remaining statements of Proposition~\ref{prop:snakes}.

The vertex lists and the lattice point tables were also checked by computer for $2\le n\le21$ (straight snakes) and $3\le n\le21$ (bent snakes).

\section{Coordinates}\label{app:coords}

\paragraph{The 27-gon of Proposition~\ref{prop:27}.} In coordinates in which $(u,v)$ denotes the point $(\sqrt3\,u/2,\,v/2)$, the vertices in order are
\[
\begin{gathered}
(1,-3),\ (-2,-6),\ (-1,-3),\ (-4,0),\ (-2,0),\ (-2,2),\ (-1,3),\ (-2,6),\ (1,3),\\
(2,6),\ (2,0),\ (6,0),\ (6,6),\ (7,3),\ (10,6),\ (9,3),\ (12,0),\ (10,0),\\
(10,-2),\ (9,-3),\ (10,-6),\ (7,-3),\ (6,-6),\ (5,-3),\ (4,-4),\ (3,-3),\ (2,-6).
\end{gathered}
\]

\paragraph{The 42-gon of Example~\ref{ex:42}.} In Cartesian coordinates, the vertices in order are
\[
\begin{gathered}
(0,2),\ (0,1),\ (-1,1),\ (-1,0),\ (-2,0),\ (-1,-1),\ (-1,-2),\ (-3,-2),\ (-3,-1),\\
(-4,-1),\ (-4,0),\ (-5,-1),\ (-5,-2),\ (-6,-2),\ (-5,-3),\ (-4,-3),\ (-4,-4),\ (-3,-3),\\
(-2,-3),\ (-2,-4),\ (-1,-3),\ (0,-4),\ (0,-3),\ (1,-3),\ (1,-2),\ (2,-2),\ (1,-1),\\
(1,0),\ (3,0),\ (3,-1),\ (4,-1),\ (4,-2),\ (5,-1),\ (5,0),\ (6,0),\ (5,1),\\
(4,1),\ (4,2),\ (3,1),\ (2,1),\ (2,2),\ (1,1).
\end{gathered}
\]

\end{document}